\documentclass[3p,11pt]{elsarticle}
\usepackage{amsmath}
\usepackage{graphicx}
\usepackage{caption} 
\usepackage{wrapfig}
\usepackage{color,soul}
\usepackage{enumitem}
\usepackage{subcaption}
\usepackage{comment}
\usepackage{amsmath}
\usepackage{times}
\usepackage{amssymb}
\usepackage{array}
\usepackage{yhmath}
\usepackage{latexsym}
\usepackage{amsthm}
\usepackage{amsfonts}
\usepackage{nomencl}
\usepackage[table]{xcolor}
\usepackage{lineno}

\makenomenclature
\usepackage{hyperref}
\usepackage{array,tabularx}
\usepackage{scalerel,stackengine}
\usepackage{graphics}
\usepackage[english]{babel}
\allowdisplaybreaks

\stackMath
\numberwithin{equation}{section}
\date{}
\newtheorem{theorem}{Theorem}[section]

\newtheorem{lemma}{Lemma}[section]
\newtheorem{definition}{Definition}[section]

\usepackage[utf8]{inputenc}

\begin{document}
\begin{frontmatter} 
\title{Formation of vacuum state and delta-shock in the solution of two-dimensional Riemann problem for zero pressure gas dynamics \tnoteref{mytitlenote}} 
\author[1]{Anamika Pandey}
\author[1]{T. Raja Sekhar}
\address[1]{Department of Mathematics, Indian Institute of Technology Kharagpur, Kharagpur, West Bengal, India} 
\cortext[mycorrespondingauthor]{Corresponding author} 
\ead{trajasekhar@maths.iitkgp.ac.in} 

\begin{abstract}
In this article, we investigate the two-dimensional pressureless Euler equations with three constant Riemann initial data. Our primary focus is on the wave interactions involving contact discontinuities and delta shocks. A distinguishing feature of the solution is the emergence of a delta shock wave which is characterized by a Dirac delta function appearing in both the density and internal energy variables. By exploiting generalized characteristic analysis, nine topologically distinct solution patterns are derived. Some of these configurations exhibit features similar to Mach-reflection and in certain cases, vacuum regions may also develop. To validate the theoretical results, numerical simulations are carried out using a semidiscrete central upwind scheme. The comparison between analytical and numerical results demonstrates excellent agreement, providing deeper insights into the complex dynamics of wave interactions in the pressureless Euler system.  \nocite{*}
\end{abstract}

\end{frontmatter}

\section{Introduction}

The mathematical modeling of various natural phenomena heavily relies on the use of hyperbolic conservation laws, particularly those that are nonlinear in nature. Nonlinear hyperbolic conservation laws play a pivotal role in the mathematical modeling of diverse natural phenomena. Among these, the Euler equations stand as a fundamental system describing the conservation principles of mass, momentum and energy. Numerous theoretical and numerical investigations have been carried out to explore the Euler equations, see \cite{Solutionpositivescheme, yzhangconjecture, Kurganovnumerical, rahulbarthwaleuler, mzafardegenerate, rahulbarthwaljmaa}. However, due to the inherent complexity of the full Euler system, a comprehensive analytical treatment remains a formidable challenge. To address these difficulties, researchers frequently investigate reduced or simplified models that retain essential features of the original system while offering greater analytical tractability. A particularly notable example is the zero-pressure gas dynamics system, also known as the pressureless Euler system or transport equations. This model captures the essential features of pressureless flows and serves as an effective approximation in several physical contexts.
\par The widely studied two-dimensional pressureless Euler equations are formally expressed as follows:
\begin{equation}\label{1}
\begin{split}
    &\frac{\partial \rho}{\partial t}+\frac{\partial(\rho u) }{\partial x}+\frac{\partial (\rho v)}{\partial y}=0,\\
    &\frac{\partial (\rho u)}{\partial t}+\frac{\partial(\rho u^2) }{\partial x}+\frac{\partial (\rho uv)}{\partial y}=0,\\
    &\frac{\partial (\rho v)}{\partial t}+\frac{\partial(\rho uv) }{\partial x}+\frac{\partial (\rho v^2)}{\partial y}=0, \; (t,x,y) \in \mathbb{R}^+ \times \mathbb{R}^2,
\end{split}
\end{equation}
where the variables $u, v$ and $\rho \geq 0$, respectively, represent x- and y-velocity components and the density. This system (\ref{1}) arises in kinetic theory from the Boltzmann equation \cite{Bouchut} and through flux-splitting techniques implemented to the compressible Euler system \cite{li1985second,agarwal1994frontiers}. System (\ref{1}) effectively describes the dynamics of free particles that stick upon collision, offering valuable insights into the development of large-scale cosmic structures \cite{Brenier, Shandarin}. Additionally, this model has found applications in simulating dusty media, where pressure is negligible \cite{kraiko1979discontinuity}. Extensive research has been conducted on the system (\ref{1}), see \cite{zwang2dzeropressure, klyushnev2023model, choudhury2014spherically, li2001delta}.
\par System (\ref{1}) serves as a pivotal framework for exploring the delta shock wave emergence and the development of vacuum states, and has consequently attracted extensive attention. Bouchut \cite{Bouchut} derived an explicit formulation for the 1-D Riemann solution with Radon measure-valued initial data and verified that in the sense of measure, this solution satisfies system (\ref{1}). Using the variational principle, Weinan et al. \cite{weinan1996generalized} derived a global weak solution for arbitrary initial data. Sheng and Zhang \cite{ShengMemoirs} employed the method of characteristics alongside the vanishing viscosity technique to construct solutions for system (\ref{1}) in the case of one dimension, and further extended their analysis in the case of two dimensions with a four-quadrant Riemann initial condition consisting solely of contact discontinuities. For the system (\ref{1}) Li and Zhang \cite{Generalizedlizhang} developed a generalized Rankine–Hugoniot condition for two-dimensional delta shocks and utilized it to establish the complete Riemann solutions. Furthermore, Li and Yang \cite{liyangDeltashocks} extended the study to $n(n\geq 3)$-dimensional configurations, characterized by two piece of constant states.
\par A delta shock represents a non-classical, nonlinear wave characterized by the presence of a Dirac delta distribution in at least one of the state variables. Physically, such structures are used to represent concentration phenomena in fluid dynamics. The concept was first identified by Tan and Zhang \cite{TanzhangfourJcase} during their investigation of a reduced form of the 2-D gas dynamics system while studying the Riemann problem. With this new notion, Tan et al. \cite{tan1994delta} subsequently proved the existence and uniqueness of delta shock under specific viscous perturbations in the 1-D framework. The development of delta shocks is closely associated with phenomena such as cavitation, concentration and vacuum formation in weak solutions. These phenomena have been rigorously validated in the context of isentropic fluids as the pressure tends to zero \cite{chen2003formation, Msunconcentrationcavitation}, and analogous results have been demonstrated for the non-isentropic full compressible Euler system \cite{ChenConcentration}. For further information related to delta shocks, the reader is referred to \cite{rahulbarthwaldeltashock, sheng1999riemann, PANDEY2025129378, Mzafardeltashock, deltashockcshen} and the references therein.
\par Another prominent feature of the system (\ref{1}) is the emergence of Mach-reflection like configurations within its solution structure, as demonstrated in \cite{jbook} through the investigation of two-dimensional Riemann problems involving four constant initial states. Cheng et al. \cite{Chengzeropressuregas} reported similar configurations arising from initial data consisting of three constant states. The Mach-reflection like patterns represent complex global wave interactions that elude capture by sequentially solving local Riemann problems at the interaction points of external waves. Instead, a holistic, global framework is necessary to ensure the consistent alignment of all exterior wave structures. Unlike classical Mach-reflection phenomena \cite{courant1999supersonic, 2001yzhengsystems} and the Guckenheimer structure observed in 2-D scalar conservation laws \cite{ZhangGeneralized, Shengscalarconservation}, these configurations arise solely from the global interaction of nonlinear waves in two dimensions, dynamics that are entirely absent in one-dimensional systems.
\par A common characteristic observed across the aforementioned multi-dimensional delta shock waves is that the Dirac delta distribution appears in only a single state variable. However, it is well established that even in media where pressure effects are negligible, energy transport remains a critical physical process. This necessitates the incorporation of energy conservation into the zero-pressure gas dynamics framework. Motivated by these considerations, we turn our attention to the system of conservation laws governing mass, momentum, and energy in two-dimensional zero-pressure gas dynamics, which is formulated as follows:
\begin{equation}\label{2}
\begin{split}
    &\frac{\partial \rho}{\partial t}+\frac{\partial(\rho u) }{\partial x}+\frac{\partial (\rho v)}{\partial y}=0,\\
    &\frac{\partial (\rho u)}{\partial t}+\frac{\partial(\rho u^2) }{\partial x}+\frac{\partial (\rho uv)}{\partial y}=0,\\
    &\frac{\partial (\rho v)}{\partial t}+\frac{\partial(\rho uv) }{\partial x}+\frac{\partial (\rho v^2)}{\partial y}=0,\\
    &\frac{\partial }{\partial t} \left(\rho e +\frac{\rho}{2}(u^2+v^2) \right)+\frac{\partial}{\partial x}\left(\left(\rho e +\frac{\rho}{2}(u^2+v^2)\right)u \right)+\frac{\partial}{\partial y} \left(\left(\rho e +\frac{\rho}{2}(u^2+v^2)\right)v \right)=0,
\end{split}
\end{equation}
where the variable $e \geq 0$ denotes specific energy. System (\ref{2}) arises by setting the pressure $p =0$ in the 2-D compressible Euler system. Nilsson et al. \cite{Massmomentumenergy, Massmomentumenergy1} analyzed the system (\ref{2}) by introducing a new variable $H=\rho e$ to represent the internal energy and established the simultaneous concentration of mass and energy along the one-dimensional delta shock front. Building on these insights, Cheng \cite{chengzeropressuregas1} provided a complete solution to the one-dimensional Riemann problem. Inspired by these contributions, Pang \cite{ypangzeropressure} studied (\ref{2}) in the following form:

\begin{equation}\label{3}
\begin{split}
    &\frac{\partial \rho}{\partial t}+\frac{\partial(\rho u) }{\partial x}+\frac{\partial (\rho v)}{\partial y}=0,\\
    &\frac{\partial (\rho u)}{\partial t}+\frac{\partial(\rho u^2) }{\partial x}+\frac{\partial (\rho uv)}{\partial y}=0,\\
    &\frac{\partial (\rho v)}{\partial t}+\frac{\partial(\rho uv) }{\partial x}+\frac{\partial (\rho v^2)}{\partial y}=0,\\
    &\frac{\partial }{\partial t} \left(H +\frac{\rho}{2}(u^2+v^2) \right)+\frac{\partial}{\partial x}\left(\left(H +\frac{\rho}{2}(u^2+v^2)\right)u \right)+\frac{\partial}{\partial y} \left(\left(H +\frac{\rho}{2}(u^2+v^2)\right)v \right)=0,
\end{split}
\end{equation}
where $H \geq 0$ denotes the internal energy and by considering four-quadrant Riemann initial conditions, he examined various solution structures along with their corresponding admissibility criteria, all within the framework of assumption (H) \cite{tan1994delta}. 
\par Based on these foundations, this article focuses on construction of the Riemann solution for the 2-D system (\ref{3}) with three constant initial states given by: 
\begin{equation}\label{4}
(\rho, u,v,H)(x,y,0)= \begin{cases}
      (\rho_1, u_1, v_1, H_1), & \text{if} \;x>0,y>0, \\
      (\rho_2, u_2, v_2, H_2), & \text{if} \;x<0,y>0,\\
      (\rho_3, u_3, v_3, H_3), & \text{if} \;y<0,x\in \mathbb{R}.
   \end{cases}
\end{equation}
In contrast to the four-state configuration studied in \cite{ypangzeropressure}, our aim is to conduct a detailed investigation of the complex interactions among contact discontinuities and delta shocks arising from a three-state initial configuration. By systematically classifying the initial data (\ref{4}) according to the number and spatial arrangement of delta shock waves, we identify three primary cases: configurations involving three, two, or a single delta shock. Employing the generalized characteristic analysis method, we construct the corresponding Riemann solutions for each case, ultimately revealing nine topologically distinct patterns of wave interactions where vacuum arises in certain cases. Notably, the configuration consisting solely of contact discontinuities yields a trivial solution. 
\par To ensure a comprehensive analysis of multidimensional Riemann problems for system \eqref{3}, it is essential to explore a wide range of initial configurations. This motivation led us to consider three-quadrant Riemann data, which balances analytical tractability with the ability to capture the system’s essential wave interaction phenomena. Despite the reduction in the number of initial states considered, our framework preserves the essential mathematical structure and intricate wave dynamics presented in \cite{ypangzeropressure}. In fact, our formulation streamlines the classification process while retaining all core phenomena, including Mach-reflection-like configurations, thereby achieving a balance between analytical clarity and the inherent complexity of wave interactions. The adoption of three-quadrant initial data significantly simplifies the analysis, reducing the 23 potential interaction patterns to just nine topologically distinct cases without compromising generality. In contrast to the purely analytical treatment of the four-state Riemann problem discussed in \cite{ypangzeropressure}, our study advances the analysis by incorporating numerical simulations through a second-order semidiscrete central-upwind scheme. This approach enables the capture of intricate features such as vacuum formation and Mach-reflection-like structures, which have not been previously addressed. To substantiate the validity of our analytical solutions, we employ the second-order semidiscrete central-upwind scheme to numerically solve the Riemann problem (\ref{3}–\ref{4}). The strong agreement between the analytical and numerical results highlights the robustness and accuracy of our constructed solutions, successfully reflecting the underlying geometric wave structures.
\par The structure of the article is as follows: Section 2 provides a detailed discussion of the essential properties of the system (\ref{3}–\ref{4}). The qualitative behavior of its solutions is analyzed in Section 3. Section 4 presents the numerical scheme adopted to verify the analytical findings and we classify the Riemann solutions based on distinct combinations of exterior waves in Section 5. These solutions are then constructed using the generalized characteristic method in Section 6. Section 7 closes the article by summarizing the central results and future works.

\section{Preliminary Analysis}
In this section, we explore the core properties of the system (\ref{3}) subject to the initial condition (\ref{4}). For a more comprehensive discussion of these foundational concepts, readers may refer to \cite{jbook, ypangzeropressure, TanzhangfourJcase}.

\subsection{\bf{Hyperbolicity}}

For smooth solutions, system (\ref{3}) can be expressed in a matrix form as follows:
\begin{equation}\label{hyp1}
\begin{aligned}
   \begin{bmatrix}
    \rho \\ u \\v \\ H
    \end{bmatrix}_t+
    \begin{bmatrix}
     u & \rho & 0 & 0 \\
     0 & u & 0 & 0\\
     0 & 0 & u & 0\\
     0 & H & 0 & u
     \end{bmatrix} \begin{bmatrix}
    \rho \\ u \\v \\ H
    \end{bmatrix}_x+\begin{bmatrix}
     v & 0 & \rho & 0 \\
     0 & v & 0 & 0\\
     0 & 0 & v & 0\\
     0 & 0 & H & v
     \end{bmatrix} \begin{bmatrix}
      \rho \\ u \\v \\ H
    \end{bmatrix}_y  & =0. 
\end{aligned}
\end{equation}
For any given direction $(\mu, \nu)$ satisfying $\mu^2+\nu^2=1$, the characteristic equation corresponding to system (\ref{hyp1}) is given by:
\begin{equation*}
\begin{aligned}
 \begin{vmatrix}
    (\mu u+\nu v)-\lambda & \mu \rho & \nu \rho & 0 \\
     0 & (\mu u+\nu v)-\lambda & 0 & 0\\
     0 & 0 & (\mu u+\nu v)-\lambda & 0\\
     0 & \mu H & \nu H & (\mu u+\nu v)-\lambda
    \end{vmatrix} & =0,
\end{aligned}
\end{equation*}
which implies that the system (\ref{3}) has four repeated eigenvalues as 
 $\begin{displaystyle}{\lambda_i=\mu u + \nu v\; ; \; i=1,2,3,4} \end{displaystyle}$ 
with the corresponding eigenvectors ${\bf{r_1}}=(1,0,0,0)^T, {\bf{r_2}}=(0,\nu,-\mu,0)^T\;\text{and}\;{\bf{r_3}}=(0,0,0,1)^T$, respectively. Hence, system (\ref{3}) is non strictly weakly hyperbolic. Further, it is evident that $\nabla \lambda_i \cdot \mathbf{r}_i = 0$ for $i = 1, 2, 3$, which follows that all $\lambda_i$'s correspond to linearly degenerate characteristic fields.

\subsection{\textbf{Pseudo-characteristics and self-similar solutions in self-similar plane}}
The system (\ref{3}) along with the initial condition (\ref{4}) remain invariant under the transformation \((x, y, t) \rightarrow (\alpha x, \alpha y, \alpha t)\) for any positive scalar $\alpha > 0$. This invariance motivates the search for self-similar solutions of the form $(\rho, u, v, H)(x, y, t) = (\rho(\xi, \eta), u(\xi, \eta), v(\xi, \eta), H(\xi, \eta))$, where the similarity variables are defined as \((\xi, \eta) = (x/t, y/t)\). 
In the case of smooth solutions, system (\ref{3}) is expressed in the following matrix form:
\begin{equation}\label{pre4}
\begin{bmatrix}
u-\xi  & \rho & 0 & 0 \\
0 & \rho(u-\xi) & 0 & 0\\
0 & 0 & \rho(u-\xi) & 0\\
0 & H & 0 & u-\xi
\end{bmatrix}\begin{bmatrix}
\rho \\ u \\ v \\ H
\end{bmatrix}_{\xi}+\begin{bmatrix}
v-\eta  & 0 & \rho & 0 \\
0 & \rho(v-\eta) & 0 & 0\\
0 & 0 & \rho(v-\eta) & 0\\
0 & 0 & H & v-\eta
\end{bmatrix}\begin{bmatrix}
\rho \\ u \\ v \\ H
\end{bmatrix}_{\eta}=0,
\end{equation}
and the initial condition (\ref{4}) transforms into:
\begin{equation}\label{pre3}
\begin{aligned}
{\lim_{\substack{\xi^{2}+\eta ^{2} \rightarrow \infty \\ {\xi / \eta =const.}}}} \left(  \rho,u,v, H\right) \left( \xi,\eta \right) =\begin{cases}(\rho_1, u_1, v_1, H_1) & ,\xi,\eta >0,\\
(\rho_2, u_2, v_2, H_2) &,\xi <0,\eta >0,\\ (\rho_3, u_3, v_3, H_3) &, \xi\in \mathbb{R} ,\eta <0.\end{cases}
\end{aligned}
\end{equation}
This system admits four identical eigenvalues given by $\Lambda_i=\dfrac{\eta - v}{\xi -u}, i=1,2,3,4$ with the corresponding eigenvectors ${\bf{r_1}}=(1,0,0,0)^T, {\bf{r_2}}=(0,u-\xi,v-\eta,0)^T$ and ${\bf{r_3}}=(0,0,0,1)^T$. It is evident that $\nabla \Lambda \cdot {\bf{r_i}}=0, i = 1,2,3$. This indicates that the $\Lambda_i$ characteristic fields are linearly degenerate. 

In $(\xi,\eta)$-plane, the pseudo-characteristic curve corresponding to system (\ref{3}) is
\begin{equation}
    \frac{d\eta}{d\xi} = \Lambda,
\end{equation}
having the singularity occurring at the point $(u,v)$. Along these curves, the velocities satisfy, the following relation
\begin{align*}
\begin{cases}
     \frac{du}{d\xi} = 0,\\ \frac{dv}{d\xi} = 0.
\end{cases}
\end{align*}
This implies that both $u$ and $v$ remain constant and therefore, all pseudo-characteristics are straight lines. We orient the integral curves of $\dfrac{d \eta}{d \xi}=\Lambda_i$ towards singularity point.

In addition to a constant state $(\rho, u, v, H)(\xi, \eta) = \text{Const.}$ with $\rho > 0$ and $H > 0$, the system (\ref{3}) also admits a vacuum type solution of the form 
\begin{align}\label{pre5}
    \text{Vac:} \;\; (\rho,u,v,H)(\xi,\eta) = (0, u(\xi,\eta), v(\xi,\eta), 0),
\end{align}
where $u(\xi, \eta)$ and $v(\xi, \eta)$ are arbitrary smooth vector-valued functions.

\subsection{\bf{Bounded discontinuity solutions}}

Let $\eta=\eta(\xi)$, corresponding to the surface $y=t\eta\left(\frac{x}{t}\right)$, be a smooth discontinuity curve associated with a bounded discontinuous solution. Assume the limit states on either side of this surface are $(\rho_R, u_R, v_R, H_R)$ on the right and $(\rho_L, u_L, v_L, H_L)$ on the left. Then the Rankine-Hugoniot conditions must be satisfied, \cite{jbook}, i.e.,

\begin{equation}\label{pre6}
\begin{cases}
(\eta-\xi\sigma)[\rho]+\sigma[\rho u]-[\rho v]=0,\\
(\eta-\xi\sigma)[\rho u]+\sigma[\rho u^2]-[\rho uv]=0,\\
(\eta-\xi\sigma)[\rho v]+\sigma[\rho uv]-[\rho v^2]=0,\\
(\eta-\xi\sigma)\left[\frac{\rho}{2}(u^2+v^2)+H \right]+\sigma \left[u\left(\frac{\rho}{2}(u^2+v^2)+H\right)\right]-\left[v\left(\frac{\rho}{2}(u^2+v^2)+H\right)\right]=0,
\end{cases}
\end{equation}
where $[G] = G_R-G_L$ denotes the jump across the discontinuity and $(\eta-\xi\sigma,\sigma,-1)$ where $\sigma=\eta'(\xi)$ is the normal to this surface.

Solving the equation (\ref{pre6}), we obtain 
\begin{align}
    J(\xi)=\begin{cases}
        \displaystyle{\sigma=\frac{d\eta}{d\xi}=\frac{\eta-v_L}{\xi-u_L}=\frac{\eta-v_R}{\xi-u_R}},\\
       \displaystyle{ [v]=\sigma[u]}.
    \end{cases}
\end{align}
We orient the integral curve of $\displaystyle{\frac{d\eta}{d\xi}}=\sigma$ towards the singularity point $\Xi_L=(u_L,v_L)$ or $\Xi_R= (u_R,v_R)$.

\subsection{\bf{Delta shock solution}}
\begin{definition}
    The three-dimensional weighted $\delta$- measure $w(t, s)\delta_s$ supported on a smooth surface parametrized as $\varPhi: ~x=x(t, s),~ y=y(t, s)$ can be defined as
\begin{align}\label{pre7}
\langle w(t, s)\delta, \phi \rangle=\displaystyle{\int_0^\infty}\displaystyle{\int_0^\infty}w(t, s)\phi(t, x(t, s), y(t, s)) ds dt,    
\end{align}
for any test function $\phi(t, x, y) \in C_0^\infty ([0, \infty)\times \mathbb{R}^2).$
\end{definition}
Let $\varPhi$ be a smooth discontinuity surface on which $\rho$ and $H$ become a Dirac measure. We parameterize the surface as $\varPhi: ~X=X(t, s),$ where $X=(x,y)$ and $s\in [0, \infty)$ is the parameter, which divides the $(t, x, y)$ space into two infinite half spaces $\Omega_-$ and $\Omega_+$. Further, we assume that $\mathbf{n}$ is normal to the surface $\varPhi$ which is oriented from $\Omega_-$ to $\Omega_+$. Here we consider the delta shock wave solution of the form \cite{tan1994delta, sheng1999riemann}         
\begin{align}\label{pre8}
(\rho,u,v,H)(t, x, y)= \begin{cases}
    (\rho_L, u_L, v_L,H_L), &(t, X)\in \Omega_-,\\
    (w(t, s)\delta(X-X(t,s)), u_\delta(t,s), v_\delta(t,s), h(t, s)\delta(X-X(t,s)) &(t, X) \in \varPhi,\\
    (\rho_R, u_R, v_R,H_R),&(t, X)\in \Omega_+,
    \end{cases}
\end{align}
where $\delta$ is the standard dirac measure and  $w(t, s), h(t,s)\in C^1([0, \infty)\times [0, \infty))$ are the weights of the delta shock on both the state variables $\rho$ and $H$, respectively. 
\begin{lemma}
A delta shock wave solution expressed in (\ref{pre8}) is considered a distributional solution of system (\ref{3}) provided it fulfills the generalized Rankine-Hugoniot conditions given by
\begin{subequations}\label{pre9}
    \begin{align}
        &\dfrac{\partial x}{\partial t}= u_\delta(t,s),\\[5pt]
        &\dfrac{\partial y}{\partial t}=v_\delta(t, s),\\[5pt]
        & \frac{\partial w}{\partial t}=(n_t,n_x,n_y)\cdot([\rho], [\rho u],[\rho v]),\\[5pt]
        &\frac{\partial (w u_\delta)}{\partial t}=(n_t,n_x,n_y)\cdot([\rho u], [\rho u^2],[\rho uv]),\\[5pt]
        &\frac{\partial (w v_\delta)}{\partial t}=(n_t, n_x, n_y).([\rho v],[\rho uv],[\rho v^2]), \\[5pt]
       & \frac{\partial (\frac{1}{2}w (u_\delta^2+v_\delta^2)+h)}{\partial t}=(n_t, n_x, n_y)\cdot\left(\left[\frac{\rho}{2}(u^2+v^2)+H \right], \left[u\left(\frac{\rho}{2}(u^2+v^2)+H\right)\right],\left[v\left(\frac{\rho}{2}(u^2+v^2)+H\right)\right]\right),
    \end{align}
\end{subequations}
in which $[G]=G_+- G_-$ denotes the jump of $G$ across the discontinuity surface $\varPhi: x=x(t, s),~y=y(t, s)$. The normal vector to the surface $\varPhi$ is $(n_t, n_x, n_y)=\left(u_\delta \dfrac{\partial y}{\partial s}-v_\delta \dfrac{\partial x}{\partial s}, -\dfrac{\partial y}{\partial s}, \dfrac{\partial x}{\partial s}\right)$.
\end{lemma}

Furthermore, to ensure a unique delta shock-type solution structured as (\ref{pre8}), the entropy conditions must also be satisfied \cite{jbook}
\begin{align}\label{pre10}
    (u_R,v_R)\cdot(n_x,n_y)<(u_\delta,v_\delta)\cdot(n_x,n_y)<(u_L,v_L)\cdot(n_x,n_y),
\end{align}
which implies that no characteristic line is outgoing on either side of discontinuity.

We follow the approach inspired by \cite{jbook} and introduce a pseudo-self-similar transformation:
\begin{align}\label{pre11}
    \begin{cases}
        \begin{split}
        X(t,s)&=t\Xi(\bar{s}),\\
        U_\delta(t,s)&=U_\delta(\bar{s}),\\
        w(t,s) &=tm(\bar{s}),\\
        h(t,s)&=tn(\bar{s}),
    \end{split}
    \end{cases}
\end{align}
where $\bar{s}=\ln(s/t)\leq 0, \Xi=(\xi,\eta)$ and $U=(u,v)$ then we can rewrite the equations (\ref{pre9}) and (\ref{pre10}) as follows
 \begin{align}\label{pre12}
        \begin{split}
            \Xi-\Xi'&=U_\delta,\\
            m-m' &=\left(\rho_L[U_L,U_\delta,\Xi]-\rho_R[U_R,U_\delta,\Xi] \right)e^{-\bar{s}},\\
            mU_\delta-(mU_\delta)' &=\left(\rho_LU_L[U_L,U_\delta,\Xi]-\rho_RU_R[U_R,U_\delta,\Xi] \right)e^{-\bar{s}},\\
            \frac{m}{2}||U_\delta^2||+n-\left(\frac{m}{2}||U_\delta^2||+n\right)' &=\left(\left(\frac{\rho_L}{2}||U_L^2||+H_L\right)[U_L,U_\delta,\Xi] -\left(\frac{\rho_R}{2}||U_R^2||+H_R\right)[U_R,U_\delta,\Xi]\right) e^{-\bar{s}},
        \end{split}
    \end{align}
and
\begin{align}\label{pre13}
    [U_L,U_\delta,\Xi]>0, \; \; \; [U_R,U_\delta,\Xi]<0,
\end{align}
where we denote
\begin{align}\label{pre14}
    [U_L,U_R,\frac{\partial X}{\partial s}]=\begin{vmatrix}
        u_L &v_L &1\\
        u_R&v_R &1\\
        \frac{\partial x}{\partial s}&\frac{\partial y}{\partial s}&0
    \end{vmatrix}, [U_L,U_R,\Xi]=\begin{vmatrix}
        u_L&v_L&1\\
        u_R&v_R&1\\
        \xi&\eta&1
    \end{vmatrix}.
\end{align}


\section{The qualitative behavior of the solution}
\begin{figure}
\begin{subfigure}{.5\textwidth}
  \centering
  \includegraphics[width=.6\linewidth]{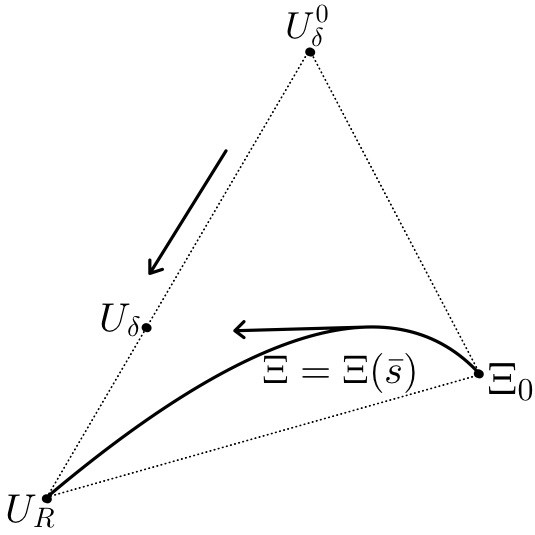}
  \caption{}  
  \label{qfig1}
\end{subfigure}%
\begin{subfigure}{.5\textwidth}
  \centering
  \includegraphics[width=.6\linewidth]{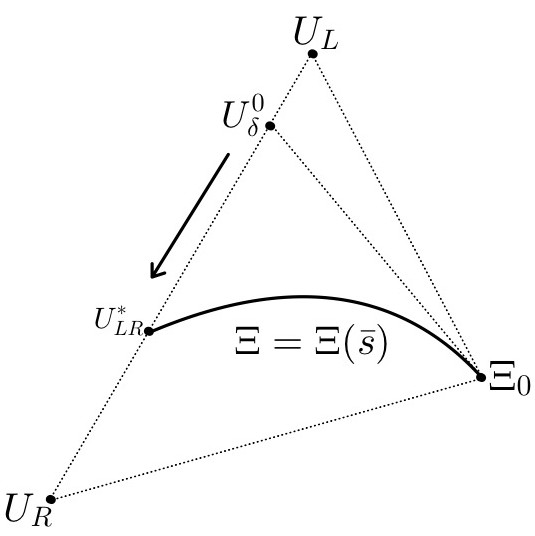}
  \caption{}  
  \label{qfig2}
\end{subfigure}%

\begin{subfigure}{.5\textwidth}
  \centering
  \includegraphics[width=.6\linewidth]{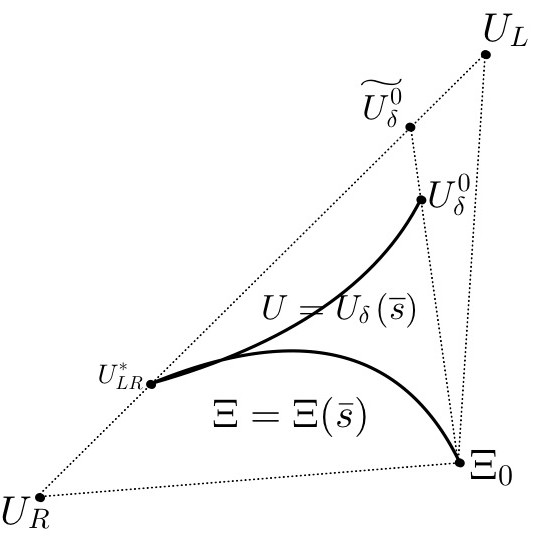}
  \caption{}  
  \label{qfig3}
\end{subfigure}%
\begin{subfigure}{.5\textwidth}
  \centering
  \includegraphics[width=.6\linewidth]{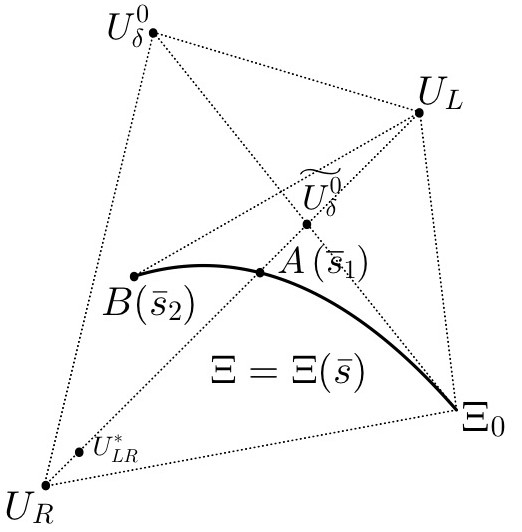}
  \caption{} 
  \label{qfig4}
\end{subfigure}
\caption{The qualitative behaviour of the solution to the reformulated generalized Rankine–Hugoniot conditions (\ref{pre12}) under the initial data (\ref{qua1})}
\end{figure}

We examine the qualitative behavior of the solution to relation (\ref{pre12}) under the given initial data
\begin{align}\label{qua1}
    \bar{s}=0: \; \Xi(0)=\Xi_0, U_\delta(0)=U_\delta^0, m(0)=m_0>0, n(0)=n_0>0,
\end{align}
which satisfies the entropy condition
\begin{align}
    [U_L,U_\delta^0,\Xi_0]>0, \; \; \; [U_R,U_\delta^0,\Xi_0]<0.
\end{align}
Hereafter, we refer to $\Xi=\Xi(\bar{s})$ as the delta shock curve, $U_\delta=U_\delta(\bar{s})$ as the velocity curve and $U^*_{LR}=\frac{U_L\sqrt{\rho_L}+U_R\sqrt{\rho_R}}{\sqrt{\rho_L}+\sqrt{\rho_R}}$ as the weighted average of the velocities $U_L$ and $U_R$. The qualitative behavior of solution to the reformulated generalized Rankine–Hugoniot relations is characterized through the following four distinct cases, following the approach adopted in \cite{ypangzeropressure}.
\begin{lemma} \label{Casea}
    When a vacuum appears as the initial state on one side of the delta shock, i.e., $\rho_RH_R>0, \rho_L=0, H_L=0$. The key features of the solution are illustrated in Figure \ref{qfig1}
\begin{itemize}
    \item $\displaystyle{\lim_{\bar{s} \to -\infty}(m,n,U_\delta,\Xi)=(0,0,U_R,U_R)}$.
    \item The $U_\delta=U_\delta(\bar{s})$ monotonously changes and approaches $U_R$ asymptotically.
    \item The convexity of $\Xi = \Xi(\bar{s})$ is clear due to its tendency to protrude to the straight line $\Xi_0U_\delta^0$.
\end{itemize}
\end{lemma}

\begin{lemma}\label{Caseb}
    When $U_\delta^0$ lies on the straight line through $U_L$ and $U_R$ i.e. $\rho_LH_L>0, \rho_RH_R>0, [U_L,U_R,U_\delta^0]=0$. The key features of the solution are demonstrated in Figure \ref{qfig2}
\begin{itemize}
    \item $\displaystyle{\lim_{\bar{s} \to -\infty}(m,n,U_\delta,\Xi)=(\sqrt{\rho_L\rho_R}[U_L,U_R,\Xi_0],n_{LR},U_{LR}^*,U_{LR}^*)}$, where $n_{LR}=\frac{[U_L,U_R,\Xi_0]}{\sqrt{\rho_L}+\sqrt{\rho_R}}(\sqrt{\rho_L}(\frac{\rho_R}{2}||U_R||^2+H_R-\frac{\rho_R}{2}||U_{LR}^*||^2)+\sqrt{\rho_R}(\frac{\rho_L}{2}||U_L||^2+H_L-\frac{\rho_L}{2}||U_{LR}^*||^2)  )$.
    \item The velocity $U_\delta=U_\delta(\bar{s})$ lies on the line $U_LU_R$ and approaches $U_{LR}^*$ asymptotically.
    \item The convexity of $\Xi = \Xi(\bar{s})$ is evident since it always protrudes to the straight line $\Xi_0U_\delta^0$.
\end{itemize}
\end{lemma}

\begin{lemma}\label{Casec}
    When $U_\delta^0$ and $\Xi_0$ are situated on the same side of the straight line joining $U_L$ and $U_R$ i.e. $\rho_LH_L>0, \rho_RH_R>0, [U_L,U_R,U_\delta^0]>0$. Also, it holds $[U_L,U_R,\Xi_0]\geq [U_L,U_R,U_\delta^0]>0$. The properties of the solution are expressed in Figure \ref{qfig3}
\begin{itemize}
    \item $\displaystyle{\lim_{\bar{s} \to -\infty}\Xi(\bar{s})=\lim_{\bar{s} \to -\infty}U_\delta(\bar{s})=U_{LR}^*}$.
    \item The $U_\delta=U_\delta(\bar{s})$ approaches $U_{LR}^*$ asymptotically and exhibits convexity opposite to that of $\Xi=\Xi(\bar{s})$.
    \item The convexity of $\Xi = \Xi(\bar{s})$ can be seen from how it protrudes to the straight line $\Xi_0U_\delta^0$.
\end{itemize}
\end{lemma}
\begin{lemma}\label{Cased}
    When both $\Xi_0$ and $U_\delta^0$ lie on the opposite sides of the straight line connecting $U_L$ and $U_R$ i.e. $\rho_LH_L>0, \rho_RH_R>0, [U_L,U_R,U_\delta^0]<0$. The behavior of the solution is demonstrated in Figure \ref{qfig4}
\begin{itemize}
    \item The curve $\Xi = E(\bar{s})$ must intersect the straight line segment joining $U_L$ and $U_R$ at a point $A$ corresponding to $\bar{s} = \bar{s}_1$. At another point $B$ with $\bar{s} = \bar{s}_2$, the entropy condition no longer holds.
    \item The segments $\widehat{\Xi_0A}$ and $\widehat{AB}$ share the same convexity which arises from the fact that they consistently protrude toward the straight line $\Xi_0U_\delta^0$. 
\end{itemize}
\end{lemma}
 

\section{Numerical discretization}
 
The central scheme provides a straightforward and flexible method for approximating solutions to nonlinear hyperbolic conservation laws. In what follows, we provide a brief overview of the semidiscrete central upwind (CU) scheme applied to the two-dimensional pressureless gas dynamics system (\ref{1}), which is expressed in conservative form:
\begin{equation}\label{n1}
\mathbf{U}_t+\mathbf{F(U)}_x+\mathbf{G(U)}_y=0, 
\end{equation}
where, 
\begin{equation}\label{n2}
\mathbf{U}=\begin{bmatrix}
    \rho\\\rho u\\ \rho v 
\end{bmatrix}, \mathbf{F(U)}=\begin{bmatrix}
    \rho u\\ \rho u^2 \\ \rho uv 
\end{bmatrix}, \mathbf{G(U)}=\begin{bmatrix}
     \rho v\\ \rho uv \\ \rho v^2  
\end{bmatrix}.
\end{equation}
Note that we do not include the energy equation in our numerical simulation since it does not affect the first three equations \cite{MR4538785}. In two dimensions, a second-order semidiscrete CU scheme is derived as follows in conservative form:
\begin{align}\label{n3}
\begin{split}
    &\frac{d}{dt}\bar{U}_{j,k}(t) = -\frac{H^x_{j+\frac{1}{2},k}(t)-H^x_{j-\frac{1}{2},k}(t)}{\Delta x}-\frac{H^y_{j,k+\frac{1}{2}}(t)-H^y_{j,k-\frac{1}{2}}(t)}{\Delta y},
    \end{split}
\end{align}
where the numerical fluxes $H^x_{j+\frac{1}{2},k}$ and $H^y_{j,k+\frac{1}{2}}$ are given by \cite{centralscheme, kurganov2001semidiscrete}
\begin{align}\label{n4}
   \begin{split}
      &  H^x_{j+\frac{1}{2},k}= \frac{a^+_{j+\frac{1}{2},k}F(U^E_{j,k})-a^-_{j+\frac{1}{2},k}F(U^W_{j+1,k})}{a^+_{j+\frac{1}{2},k}-a^-_{j+\frac{1}{2},k}}+a^+_{j+\frac{1}{2},k}a^-_{j+\frac{1}{2},k}\left[ \frac{U^W_{j+1,k}-U^E_{j,k}}{a^+_{j+\frac{1}{2},k}-a^-_{j+\frac{1}{2},k}}-q^x_{j+\frac{1}{2},k}\right],\\
     \text{and} &\\
     &H^y_{j,k+\frac{1}{2}}= \frac{b^+_{j,k+\frac{1}{2}}G(U^N_{j,k})-b^-_{j,k+\frac{1}{2}}G(U^S_{j,k+1})}{b^+_{j,k+\frac{1}{2}}-b^-_{j,k+\frac{1}{2}}}+b^+_{j,k+\frac{1}{2}}b^-_{j,k+\frac{1}{2}}\left[ \frac{U^S_{j,k+1}-U^N_{j,k}}{b^+_{j,k+\frac{1}{2}}-b^-_{j,k+\frac{1}{2}}}-q^y_{j,k+\frac{1}{2}}\right].
   \end{split}
\end{align}
It is important to note that all quantities in equation (\ref{n4}) are functions of time $t$; however, for simplicity, this dependence is omitted in the notation. 

In equation (\ref{n4}), the point values $U^{E(W,S,N)}$ are obtained from a conservative, non-oscillatory piecewise polynomial reconstruction of the desired order. For instance, a second-order CU scheme employs a piecewise linear reconstruction.
\begin{align}\label{n5}
    \begin{split}
        & U^{E(W)}_{j,k} = \bar{U}_{j,k}(t)\pm \frac{\Delta x}{2} (U_x)_{j,k},\\
        & U^{N(S)}_{j,k} = \bar{U}_{j,k}(t)\pm \frac{\Delta y}{2} (U_y)_{j,k}.
    \end{split}
\end{align}
To maintain the non-oscillatory behaviour of the reconstruction and consequently the second-order CU scheme, the slopes in (\ref{n5}) are determined by exploiting nonlinear minmod limiter defined as \cite{nessyahu1990non, torobook, Leveque, hesthaven}
\begin{align}\label{n6}
    \begin{split}
       & (U_x)_{j,k} = \mathrm{minmod}\left(\theta \frac{\bar{U}_{j+1,k}-\bar{U}_{j,k}}{\Delta x}, \frac{\bar{U}_{j+1,k}-\bar{U}_{j-1,k}}{2\Delta x}, \theta \frac{\bar{U}_{j,k}-\bar{U}_{j-1,k}}{\Delta x} \right),\\
        & (U_y)_{j,k} = \mathrm{minmod}\left(\theta \frac{\bar{U}_{j,k+1}-\bar{U}_{j,k}}{\Delta y}, \frac{\bar{U}_{j,k+1}-\bar{U}_{j,k-1}}{2\Delta y}, \theta \frac{\bar{U}_{j,k}-\bar{U}_{j,k-1}}{\Delta y} \right),     
    \end{split}
\end{align}
 where $\theta \in [1,2]$ and the multivariate minmod function is defined by,
\begin{align}
    \begin{split}
       & \mathrm{minmod}(a_1,a_2,a_3,\cdots)=\begin{cases}
        \min(a_1,a_2,a_3,\cdots) \; \; \;  \text{if}   \;  \; \;a_i>0 \; \forall \;i,\\
        \max(a_1,a_2,a_3,\cdots) \; \; \;  \text{if}   \;  \; \;a_i<0 \; \forall \;i,\\
        0  \;  \; \;  \; \;  \; \;  \;  \; \;  \; \;  \; \;  \; \; \;  \; \;  \; \;  \; \;  \;  \; \;  \; \;  \; \;\text{otherwise}.
    \end{cases}
    \end{split}
\end{align}
As all the eigenvalues of the Jacobian matrices $\frac{\partial F}{\partial U}$ and $\frac{\partial G}{\partial U}$ have multiplicity $4$ and are equal to $u$ and $v$, respectively. The one-sided local speeds in (\ref{n4}) can be estimated straightforwardly:
\begin{align}\label{n7}
    \begin{split}
        & a^+_{j+\frac{1}{2},k} = \max(u^W_{j+1,k},u^E_{j,k},0), \; a^-_{j+\frac{1}{2},k} = \min(u^W_{j+1,k},u^E_{j,k},0),\\
        & b^+_{j,k+\frac{1}{2}} = \max(v^S_{j,k+1},v^N_{j,k},0), \; b^-_{j,k+\frac{1}{2}} = \min(v^S_{j,k+1},v^N_{j,k},0).
    \end{split}
\end{align}
Now, we define the anti-diffusion term added in the flux evaluation (\ref{n4}) to reduce numerical dissipation as
\begin{align}
    \begin{split}
       & q^x_{j+\frac{1}{2},k} = \mathrm{minmode} \left( \frac{U^{NW}_{j+1,k}-U^{\text{int}}_{j+\frac{1}{2},k}}{a^+_{j+\frac{1}{2},k}-a^-_{j+\frac{1}{2},k}}, \frac{U^{\text{int}}_{j+\frac{1}{2},k}-U^{NE}_{j,k}}{a^+_{j+\frac{1}{2},k}-a^-_{j+\frac{1}{2},k}}, \frac{U^{SW}_{j+1,k}-U^{\text{int}}_{j+\frac{1}{2},k}}{a^+_{j+\frac{1}{2},k}-a^-_{j+\frac{1}{2},k}}, \frac{U^{\text{int}}_{j+\frac{1}{2},k}-U^{SE}_{j,k}}{a^+_{j+\frac{1}{2},k}-a^-_{j+\frac{1}{2},k}}\right),\\
        & q^y_{j,k+\frac{1}{2}} = \mathrm{minmode} \left( \frac{U^{SW}_{j,k+1}-U^{\text{int}}_{j,k+\frac{1}{2}}}{b^+_{j,k+\frac{1}{2}}-b^-_{j,k+\frac{1}{2}}}, \frac{U^{\text{int}}_{j,k+\frac{1}{2}}-U^{NW}_{j,k}}{b^+_{j,k+\frac{1}{2}}-b^-_{j,k+\frac{1}{2}}}, \frac{U^{SE}_{j,k+1}-U^{\text{int}}_{j,k+\frac{1}{2}}}{b^+_{j,k+\frac{1}{2}}-b^-_{j,k+\frac{1}{2}}}, \frac{U^{\text{int}}_{j,k+\frac{1}{2}}-U^{NE}_{j,k}}{b^+_{j,k+\frac{1}{2}}-b^-_{j,k+\frac{1}{2}}}\right),\\
        \text{where, the} \; & \text{intermediate values are}\\
        & U^{\text{int}}_{j+\frac{1}{2},k} = \frac{a^+_{j+\frac{1}{2},k}U^W_{j+1,k}-a^-_{j+\frac{1}{2},k}U^E_{j,k}-\left( F(U^W_{j+1,k})-F(U^E_{j,k})\right)}{a^+_{j+\frac{1}{2},k}-a^-_{j+\frac{1}{2},k}},\\
        & U^{\text{int}}_{j,k+\frac{1}{2}} = \frac{b^+_{j,k+\frac{1}{2}}U^S_{j,k+1}-b^-_{j,k+\frac{1}{2}}U^N_{j,k}-\left( G(U^S_{j,k+1})-G(U^N_{j,k})\right)}{b^+_{j,k+\frac{1}{2}}-b^-_{j,k+\frac{1}{2}}},
    \end{split}
\end{align}
and at the cell corners, the point values are
\begin{align}
    \begin{split}
        & U^{NE(NW)}_{j,k} = \bar{U}_{j,k}(t)\pm \frac{\Delta x}{2} (U_x)_{j,k}+\frac{\Delta y}{2} (U_y)_{j,k},\\
        & U^{SE(SW)}_{j,k} = \bar{U}_{j,k}(t)\pm \frac{\Delta x}{2} (U_x)_{j,k}-\frac{\Delta y}{2} (U_y)_{j,k}.
    \end{split}
\end{align}
For our numerical simulation, we use a uniform grid comprising $200 \times 200$ cells over the computational domain $[-0.5,0.5]\times [-0.5,0.5]$ with equal step sizes $\Delta x$ and $\Delta y$, respectively, along the $x$- and $y$-axes. The time domain $[0,T]$,where $T$ denotes the final simulation time, is discretized with a time step $\Delta t$, determined by a CFL number of $0.1$. In all simulations, we take $T=0.20$. For more details on second-order accurate numerical schemes, see, viz. \cite{cui2024optimal, 2ndorder, Hwenoscheme} and the references cited therein.

\section{Classification}
To construct the solution of the Riemann problem (\ref{3})–(\ref{4}), we classify the wave patterns into three main categories based on the number of delta shock waves appearing and their interactions: configurations involving three, two, or a single delta shock. Each category is further subdivided depending on the nature of interaction among the delta shocks, whether they intersect or evolve independently. Specifically, in the first three cases, all three delta shocks are present; in these, the difference lies in how they interact, ranging from no interaction to partial interactions. Similar reasoning applies to the remaining configurations. This leads to a total of nine distinct cases, each representing a topologically different structure of the Riemann solution.
\subsection{Cases involving three delta shock waves}
\begin{enumerate}
    \item $\delta_{12}+ \delta_{23}+\delta_{31}\; (u_1<u_2<u_3, v_2<v^{23}_{\delta}<v_1<v_3)$
    \item $\delta_{12}+ \delta_{23}+\delta_{31}\; (u_1<u_2<u_3, v_2<v_1<v^{23}_{\delta}<v_3, [\Xi_1,\Xi_3,A]>0)$, where $A$ denotes the interaction point of the delta shocks $\delta_{12}$ and $\delta_{23}$
    \item $\delta_{12}+ \delta_{23}+\delta_{31}\; (u_3<u_1<u_2, v_2<v_1<v_3, [\Xi_1,\Xi_2,\Xi_3]<0)$
\end{enumerate}
\subsection{Cases involving two delta shock waves}
\begin{enumerate}
    \item $J_{12}+ \delta_{23}+\delta_{31}\; (u_3<u_1=u_2, v_2<v_1<v_3)$
    \item $J_{12}+ \delta_{23}+\delta_{31}\; (u_1=u_2<u_3, v_2<v^{23}_{\delta}<v_1<v_3)$
    \item $J_{12}+ \delta_{23}+\delta_{31}\; (u_1=u_2<u_3, v_2<v_1<v^{23}_{\delta}<v_3)$
\end{enumerate}
\subsection{Cases involving single delta shock wave}
\begin{enumerate}
   \item $\delta_{12}+J_{23}+ J_{31}\; (u_1<u_2, v_1=v_2=v_3)$
   \item $J_{12}+ \delta_{23}+J_{31}\; (u_3<u_1=u_2, v_2<v_1=v_3)$
   \item $J_{12}+ \delta_{23}+J_{31}\; (u_1=u_2<u_3, v_2<v_1=v_3)$
\end{enumerate}


\section{Solution construction of the Riemann problem}

We proceed by employing characteristic analysis to derive solutions of (\ref{3}) and (\ref{4}) in the $(\xi, \eta)$-plane. To begin, we present Fixed point theorem of Schauder \cite{subrahmanyam2018elementary} along with some necessary notations.
\begin{theorem}\label{thm1}
    \textbf{Fixed Point Theorem of Schauder:} 
Let \( K \) be a nonempty, closed and convex subset of a Banach space \( \mathcal{B} \) and let \( T: K \to K \) be a continuous operator such that \( T(K) \) is precompact in \(\mathcal{B}\). Then, there exists a point \( x \in K \) such that \( Tx = x \).
\end{theorem}
Now, we present some notations that can be used throughout this article.
\begin{itemize}
    \setlength\itemsep{0.1em} 
    \renewcommand\labelitemi{} 
     \item $\text{\textcircled{\raisebox{-0.9pt}{i}}}$: the state $(\rho_i, u_i,v_i,H_i)$ for $i=1,2,3$
    \item \(\Xi_i \): the point \( (\xi_i, \eta_i) = (u_i,v_i)\) for $i=1,2,3$
    \item \(\Xi_i\Xi_j \): the straight line passing through \( \Xi_i \) and \( \Xi_j \)
   \item \( \overline{\Xi_i\Xi_j} \): the segment connecting the points \( \Xi_i \) and \( \Xi_j \)
    \item \( \Omega_i \):  the determination domain of the state $\text{\textcircled{\raisebox{-0.9pt}{i}}}:$
 where,
         \begin{align*}
             \begin{split}
                & \Omega_1 = \{(\xi,\eta) | \xi>u_1, \eta>v_1\},\\
                & \Omega_2 = \{(\xi,\eta) | \xi<u_2, \eta>v_2\},\\
                & \Omega_3 = \{(\xi,\eta) | \xi\in \mathbb{R}, \eta<v_3\}.
             \end{split}
         \end{align*}
    \item  \(J_{ij}\): contact discontinuity connecting states $\text{\textcircled{\raisebox{-0.9pt}{i}}}$
 and $\text{\textcircled{\raisebox{-0.9pt}{j}}}$
    \item \( \delta_{ij} \): delta shock wave joining the states $\text{\textcircled{\raisebox{-0.9pt}{i}}}$
 and $\text{\textcircled{\raisebox{-0.9pt}{j}}}$

     \item \( \delta^A_{ij} \):  delta shock wave originating from point $A$ and joining the states $\text{\textcircled{\raisebox{-0.9pt}{i}}}$ and $\text{\textcircled{\raisebox{-0.9pt}{j}}}$
     \item \( \delta^A_{i} \):  delta shock wave originating from point $A$ and joining vacuum and state $\text{\textcircled{\raisebox{-0.9pt}{i}}}$
     \item \(m_{ij}\) : the weight of delta shock $\delta_{ij}$ on $\rho$ variable
     \item \(n_{ij}\) : the weight of delta shock $\delta_{ij}$ on $H$ variable
\end{itemize}
 \subsection{Cases involving three delta shocks}
 This case occurs whenever the initial data holds $u_1<u_2, v_1<v_3, v_2<v_3$. Without loss of generality, assume $v_2<v_1$. Hence, we have $u_1<u_2, v_2<v_1<v_3$. The delta shock waves connecting the states $\text{\textcircled{\raisebox{-0.9pt}{1}}}$ and $\text{\textcircled{\raisebox{-0.9pt}{2}}}$
, $\text{\textcircled{\raisebox{-0.9pt}{2}}}$ and $\text{\textcircled{\raisebox{-0.9pt}{3}}}$
 and $\text{\textcircled{\raisebox{-0.9pt}{3}}}$ and $\text{\textcircled{\raisebox{-0.9pt}{1}}}$
 are
 \begin{align}\label{case_1}
     \begin{split}
         \delta_{12}: &\begin{cases}
             \xi = \frac{\sqrt{\rho_1}u_1+\sqrt{\rho_2}u_2}{\sqrt{\rho_1}+\sqrt{\rho_2}},
             U_\delta = \frac{\sqrt{\rho_1}U_1+\sqrt{\rho_2}U_2}{\sqrt{\rho_1}+\sqrt{\rho_2}},\\
             m_{12} = \sqrt{\rho_1\rho_2}(u_2-u_1),
             n_{12} = \left( \frac{\rho_1\rho_2||U_1-U_2||^2}{2(\sqrt{\rho_1}+\sqrt{\rho_2})} +\sqrt{\rho_2}H_1+\sqrt{\rho_1}H_2\right)\frac{u_2-u_1}{\sqrt{\rho_1}+\sqrt{\rho_2}},\\
             u_1<u_2,
         \end{cases}\\
         \delta_{23}: &\begin{cases}
             \eta = \frac{\sqrt{\rho_2}v_2+\sqrt{\rho_3}v_3}{\sqrt{\rho_2}+\sqrt{\rho_3}},
             U_\delta = \frac{\sqrt{\rho_2}U_2+\sqrt{\rho_3}U_3}{\sqrt{\rho_2}+\sqrt{\rho_3}},\\
             m_{23} = \sqrt{\rho_2\rho_3}(v_3-v_2),
             n_{23} = \left( \frac{\rho_2\rho_3||U_2-U_3||^2}{2(\sqrt{\rho_2}+\sqrt{\rho_3})} +\sqrt{\rho_2}H_3+\sqrt{\rho_3}H_2\right)\frac{v_3-v_2}{\sqrt{\rho_2}+\sqrt{\rho_3}},\\
             v_2<v_3,
             \end{cases}\\
             \text{and}\;\; \; \;\;\; \;\; \; \;\;\;\;\; \; \;\;\;\;\; \; \;\;\;\;\; \; \;\;\;&\\
        \delta_{31}: &\begin{cases}
             \eta = \frac{\sqrt{\rho_3}v_3+\sqrt{\rho_1}v_1}{\sqrt{\rho_1}+\sqrt{\rho_3}},
             U_\delta = \frac{\sqrt{\rho_1}U_1+\sqrt{\rho_3}U_3}{\sqrt{\rho_1}+\sqrt{\rho_3}},\\
             m_{31} = \sqrt{\rho_1\rho_3}(v_3-v_1),
             n_{31} = \left( \frac{\rho_1\rho_3||U_3-U_1||^2}{2(\sqrt{\rho_1}+\sqrt{\rho_3})} +\sqrt{\rho_1}H_3+\sqrt{\rho_3}H_1\right)\frac{v_3-v_1}{\sqrt{\rho_1}+\sqrt{\rho_3}},   \\
             v_1<v_3,
         \end{cases}
     \end{split}
 \end{align}
 respectively. For this case, the solution involves the following three subcases: 
 
\subsubsection{Case 1:}
For this setup the initial data satisfies $u_1<u_2<u_3, v_2<v^{23}_{\delta}<v_1$. Furthermore, suppose that $\delta_{31}$ and $\delta_{12}$ do not interact, then the three delta shock waves evolve independently and form a domain in which no pseudo-characteristics enter; as a result, a vacuum develops in that region. Denote the interaction point of $\delta_{12}$ and $\overline{\Xi_1\Xi_3}$ as $A(\xi_1,\eta_1)$, the interaction point of $\delta_{23}$ and $\overline{\Xi_1\Xi_2}$ as $B(\xi_2,\eta_2)$ and the interaction point of $\delta_{31}$ and $\overline{\Xi_2\Xi_3}$ as $C(\xi_3,\eta_3)$. Also, denote $ D=(u_\delta^{12},v_1)$ and $E=(u_2,v_\delta^{23})$. 

By solving locally the Riemann problem at the interaction point, if we analyze the solution, it becomes evident that these exterior waves cannot be globally matched.
Therefore, it is necessary to investigate how the delta shocks are matched. We now proceed with a rigorous construction of the solution. As the velocities and weights of the delta shock waves change at their points of interaction, the delta shocks are no longer able to follow their earlier routes.

We begin our analysis from the delta shock $\delta_{12}$. Let us consider a point $\Xi_0^1 \in \overline{AD}$, along with $\bar{U}_\delta^{12} \in \triangle A\Xi_1D$, and the weights $m_{12}^* \in [0, M_*]$, $n_{12}^* \in [0, N_*]$, where $M_*$ and $N_*$ are constants to be specified subsequently. We now examine the solution to equation (\ref{pre12}) subject to the following initial data:
\begin{align}\label{Case1_1}
    \bar{s}=0: \; \Xi(0)=\Xi_0^1, U_\delta(0)=U_{\delta_2}^0, m(0)=m_0^1, n(0)=n_0^1,
\end{align}
where, 
\begin{align}\label{Case1_2}
    m_0^1= m_{12}+m^*_{12},  n_0^1=n_{12}+n^*_{12},U_{\delta_2}^0= \frac{m_{12}U_\delta^{12}+m^*_{12}\bar{U}_\delta^{12}}{m_{12}+m^*_{12}}.
\end{align}
Since both mass and momentum are preserved at the point of interaction. Following standard computation, we obtain
\begin{align}\label{Case1_3} 
    \begin{split}
       \begin{cases}
            m &= \sqrt{(m_0^1)^2-2\rho_2m_0^1[U_2,U_{\delta_2}^0,\Xi_0^1](e^{-\bar{s}}-1)} \cdot e^{\bar{s}}, \\
        n & = \left(\frac{||U_{\delta_2}^0||^2}{2}m_0^1+n_0^1-\left(\frac{||U_2||^2}{2}+\frac{H_2}{\rho_2}\right)m_0^1 \right)e^{\bar{s}} +\left(\frac{||U_2||^2}{2}+\frac{H_2}{\rho_2}-\frac{||U_{\delta_2}||^2}{2}\right)m, \\
        U_\delta &= U_2+ \frac{m_0^1(U_{\delta_2}^0-U_2)e^{\bar{s}}}{m},\\
        \Xi &=U_2+e^{\bar{s}}(\Xi_0^1-U_2)-\frac{U_{\delta_2}^0-U_2}{\rho_2[U_2,U_{\delta_2}^0,\Xi_0^1]}(m-m_0^1e^{\bar{s}}). 
       \end{cases}
    \end{split}
\end{align}
Clearly, $u_{\delta_2}^0<u_\delta^{12}$. Based on Lemma \ref{Casea}, the delta shock, named $\delta_2^{\Xi_0^1}$, which separates the state $\text{\textcircled{\raisebox{-0.9pt}{2}}}$ from the vacuum, lies within the region $\triangle\Xi_0^1U_{\delta_2}^0\Xi_2$. Consequently, it must interact with $\delta_{23}$ at some point, labeled $\Xi_0^2,$ where the masses and the velocity associated with $\delta_2^{E_0^1}$ are given by $m_2^*, n_2^*$ and $U_{\delta_2}^*$, respectively. It is evident that $\Xi_0^2\in BE$ and $U_{\delta_2}^*\in  \triangle\Xi_2BE$, and we assume $m_2^*\in [0,M_*]$ and $n_2^*\in [0,N_*]$. Also, let us consider
\begin{align}\label{Case1_4}
     m_* &= \max \{ m_{12}, m_{23}, m_{31} \},  \\
     n_* &= \max \{ n_{12}, n_{23}, n_{31} \},\\
     \rho_* &= \max \{ \rho_1, \rho_2, \rho_3 \},\\
     p_* &= \max \left\{ \frac{\sqrt{\rho_3}(v_3 - v_2)}{(\sqrt{\rho_2} + \sqrt{\rho_3})(v_1 - v_2)}, \frac{\sqrt{\rho_1}}{\sqrt{\rho_1} + \sqrt{\rho_3}}, 
     \frac{\sqrt{\rho_2}(u_2 - u_1)}{(\sqrt{\rho_2} + \sqrt{\rho_1})(u^{*}_{13} - u_1)} \right\},\\
     q_*&= \max \left\{ \frac{||U_1||}{2}+\frac{H_1}{\rho_1}, \frac{||U_2||}{2}+\frac{H_2}{\rho_2}, \frac{||U_3||}{2}+\frac{H_3}{\rho_3} \right\},\\
     L_0&= \max \left\{ \max_{U^*_{\delta_2}\in  \triangle \Xi_2BE} ||U^*_{\delta_2}||, ||U_1||, ||U_2||, ||U_3||\right\}.
\end{align}

We now establish the existence of $M_*$ and $N_*$ and validate our claim. To this end, we begin by proving the following: 

Indeed, for the first case, noting that $v_\delta^{23}<v_1$ it follows that 
$$ \sqrt{\rho_3}(v_3 - v_1)< \sqrt{\rho_2}(v_1 - v_2).$$
Therefore, using this, we find that
\begin{align*}
\frac{\sqrt{\rho_3}(v_3 - v_2)}{(\sqrt{\rho_2} + \sqrt{\rho_3})(v_1 - v_2)} 
< \frac{\sqrt{\rho_3}(v_3 - v_2)}{\sqrt{\rho_3}(v_3 - v_1) + \sqrt{\rho_3}(v_1 - v_2)} = 1.
\end{align*}
For the second case, since $ \delta_{13}$ and $\delta_{12}$ do not interact, we deduce that $u^{*}_{12}<u^{*}_{13}$. It follows that 
$$(\sqrt{\rho_1}u_1 + \sqrt{\rho_2}u_2) < (\sqrt{\rho_1} + \sqrt{\rho_2})u^{*}_{13}.$$
Therefore, using this, we obtain
\begin{align*}
\frac{\sqrt{\rho_2}(u_2 - u_1)}{(\sqrt{\rho_2} + \sqrt{\rho_1})(u^{*}_{13} - u_1)} 
< \frac{\sqrt{\rho_2}(u_2 - u_1)}{\sqrt{\rho_1}u_1 + \sqrt{\rho_2}u_2 - (\sqrt{\rho_2} + \sqrt{\rho_1})u_1} = 1.
\end{align*}
Consider that $Ar$ is the measure of $\triangle \Xi_1 \Xi_2 \Xi_3$. Hence we obtain
\begin{align}\label{Case1_5}
\lim_{M \to +\infty} \left( \frac{(M + m_*)^2 p_*^2 + 2\rho_*(M + m_*)Ar}{M^2} \right) 
= p_*^2 \leq 1.
\end{align}
Hence, for a suitably large constant $M_*>0$ we have
\begin{align}\label{Case1_6}
    (M_* + m_*)^2 p_*^2 + 2\rho_*(M_* + m_*)Ar \leq M_*^2.
\end{align}
In a similar way, we have
\begin{align}\label{Case1_7}
\lim_{N \to +\infty} \left( \frac{(N + n_*) p_*+q_*M_*}{N} + \frac{L_0p^*(M_* + m_*)}{2N} \right) 
= p_* \leq 1,
\end{align}
we can find $N_*>0$ such that
\begin{align}\label{Case1_8}
    (N + n_*) p_*+q_*M_* + L_0p^*\frac{M + m_*}{2}\leq N_*.
\end{align}
We now proceed to verify that $m_* \in [0, M_*]$ and $n_* \in [0, N_*]$. For this purpose, we determine the value of $\bar{s}$, referred to as $\bar{s}_2$, which is associated with the point $\Xi_0^2$. By the equation (\ref{Case1_3}), we obtain the following expression:
\begin{align}\label{Case1_9}
v_2 + e^{\bar{s}_2}(\eta_0^1 - v_2) -\frac{v_{\delta_2}^0 - v_2}{\rho_2[U_2, U_{\delta_2}^0, \Xi_0^1]}(m(\bar{s}_2) - m_0^1 e^{\bar{s}_2}) = v_\delta^{23} = \frac{\sqrt{\rho_2}v_2 + \sqrt{\rho_3}v_3}{\sqrt{\rho_2} + \sqrt{\rho_3}}.
\end{align}
Since the derivative $(e^{-2\bar{s}} m^2(\bar{s}))'<0$, we deduce that
$$(m_0^1)^2<e^{-2\bar{s}_2} m^2(\bar{s}_2),$$
which gives
$$m_0^1 e^{\bar{s}_2}<m(\bar{s}_2).$$
Moreover, as $v_2<v_{\delta_2}^0$  it follows that
$$v_2 + e^{\bar{s}_2}(\eta_0^1 - v_2) < \frac{\sqrt{\rho_2}v_2 + \sqrt{\rho_3}v_3}{\sqrt{\rho_2} + \sqrt{\rho_3}}.$$
Hence we have
\begin{align}\label{Case1_10}
    0\leq e^{\bar{s}_2}< \frac{\sqrt{\rho_3}(v_3-v_2)}{(\sqrt{\rho_2} + \sqrt{\rho_3})(\eta_0^1-v_2)}<\frac{\sqrt{\rho_3}(v_3-v_2)}{(\sqrt{\rho_2} + \sqrt{\rho_3})(v_1-v_2)} \leq p^*\leq 1.
\end{align}
Next, observing that $-[U_2, U_{\delta_2}^0, \Xi_0^1] < Ar$, and using equations (\ref{Case1_3}), (\ref{Case1_7}) and (\ref{Case1_8}), we obtain
\begin{figure}[h!]
\centering
\begin{subfigure}{.5\textwidth}
  \centering
  \includegraphics[width=.75\linewidth]{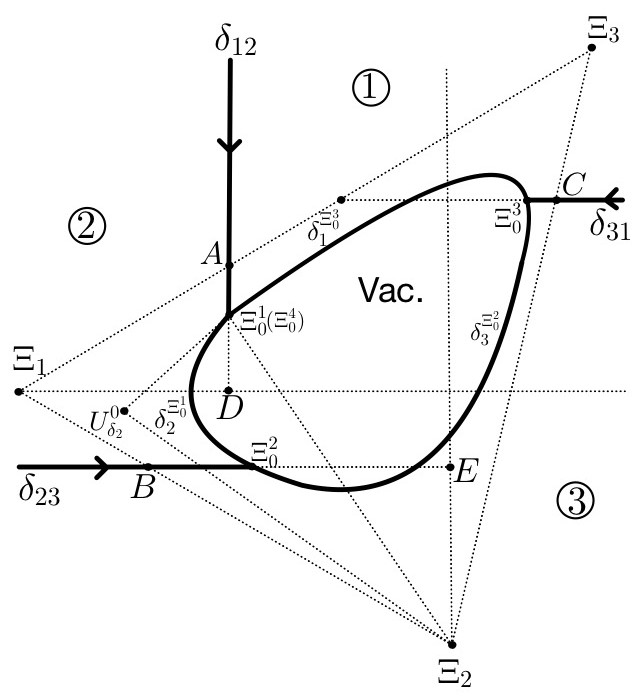}
  \caption{Analytical solution}
\end{subfigure}%
\begin{subfigure}{.5\textwidth}
  \centering
  \includegraphics[width=.7\linewidth]{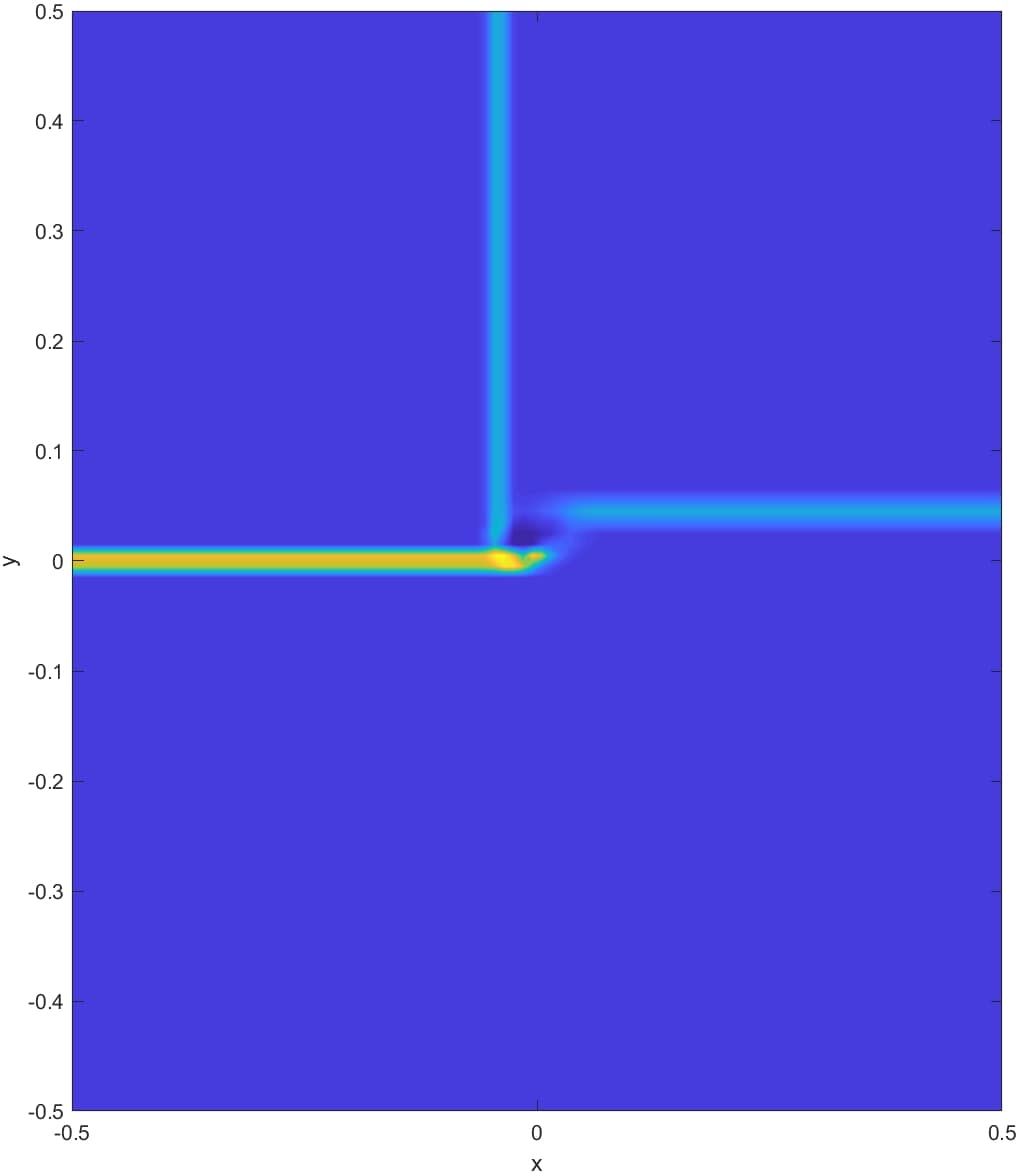}
  \caption{Numerical solution}
\end{subfigure}
\caption{Analytical and numerical solution of the Riemann problem (\ref{3})–(\ref{4}) for the initial configuration satisfying $u_1 < u_2 < u_3, v_2 < v^{23}_{\delta} < v_1<v_3$. The solution features three non-interacting delta shock waves that evolve independently, leading to the emergence of a vacuum region}
\label{case1}
\end{figure}
\begin{align}\label{Case1_11}
    \begin{split}
        m_2^*&= \sqrt{(m_0^1)^2e^{2\bar{s_2}}-2\rho_2m_0^1[U_2,U_{\delta_2}^0,\Xi_0^1](e^{\bar{s_2}}-e^{2\bar{s_2}})}\\
        & \quad \leq \sqrt{(M_* + m_*)^2 p_*^2 + 2\rho_*(M_* + m_*)Ar} \leq M_*,\\
        n_2^* &= n_0^1e^{\bar{s}_2}+\frac{||U_{\delta_2}^0||}{2}m_0^1e^{\bar{s}_2}+\left( \frac{||U_2||^2}{2}+\frac{H_2}{\rho_2}\right)(m_2^*-m_0^1e^{\bar{s}_2})-\frac{||U_{\delta_2}^*||^2}{2}m_2^* \\
        & \quad \leq (N + n_*) p_*+q_*M_* + L_0p^*\frac{M + m_*}{2}\leq N_*.
    \end{split}
\end{align}
Therefore, we establish the desired result: $0\leq m_2^*\leq M_*$ and $0 \leq n_2^* \leq N_*$.\\
Following a similar approach, we now proceed to solve equation (\ref{pre12}) with the corresponding initial data
\begin{align}\label{Case1_12}
    \bar{s}=0: \; \Xi(0)=\Xi_0^2, U_\delta(0)=U_{\delta_3}^0, m(0)=m_0^2, n(0)=n_0^2,
\end{align}
where, $$m_0^2= m_{23}+m^*_{2},  n_0^2=n_{23}+n^*_{2},U_{\delta_2}^0= \frac{m_{23}U_\delta^{23}+m^*_{2}U_{\delta_2}^*}{m_{23}+m^*_2}.$$
In a similar manner, we construct a delta shock $\delta_3^{\Xi_0^2}$ that interacts with $\delta_{31}$ at a point $\Xi_0^3$, where the associated masses and velocity of $\delta_3^{\Xi_0^2}$ are denoted by $m_3^*, n_3^*$ and $U_{\delta_3}^*$, respectively. Continuing this process, we obtain another delta shock $\delta_1^{\Xi_0^3}$ which interacts with $\delta_{12}$ at the point $\Xi_0^4 \in \overline{AD}$. The corresponding masses and velocity of $\delta_1^{\Xi_0^3}$ are $m_1^*, n_1^*$ and $U_{\delta_1}^*$, respectively. It is evident that $U_{\delta_1}^* \in \triangle A \Xi_1 D$. We further assert that $m_3^*, m_1^* \in [0, M_*]$ and $n_3^*, n_1^* \in [0, N_*]$. The proof follows similarly to that of $m_2^*, n_2^*$.

Let $K = K_1 \times K_2 \times K_3  \times K_4$, where $K_1 = AD, K_2 = [0,M_*],K_3 = [0,N_*]$ and $ K_4 = \triangle ADU_1$. Consider the operator $T: K \rightarrow K$ such that
\[
T\left(\Xi_0^1, m_{12}^*,n_{12}^*, \bar{U}_{\delta}^{12} \right) = \left(\Xi_0^4, m_1^*, n_1^*, U_{\delta_1}^* \right).
\]
We assert the existence of a fixed point $(\Xi^*, m_*,n_*, U_\delta^*) \in K$ such that $$T(\Xi^*, m_*,n_*, U_\delta^*)=(\Xi^*, m_*,n_*, U_\delta^*).$$ 

Consider $\mathbb{R}^6$ as the Banach space with the standard metric. Let $K \subset \mathbb{R}^6$ be a non-empty, convex and closed subset in $\mathbb{R}^6$. By the continuous dependence of solutions to ordinary differential equations on their initial conditions, it follows that the operator $T$ is continuous. Furthermore, as $T(K)$ is bounded in $\mathbb{R}^6$, it is precompact in $\mathbb{R}^6$. Therefore, by invoking Theorem \ref{thm1}, we conclude that $T$ admits a fixed point $(\Xi^*, U_\delta^*, m^*, n^*) \in K$ such that
\[
T(\Xi^*, U_\delta^*, m^*, n^*) = (\Xi^*, U_\delta^*, m^*, n^*).
\]
\begin{figure}
\begin{subfigure}{.5\textwidth}
  \centering
  \includegraphics[width=.75\linewidth]{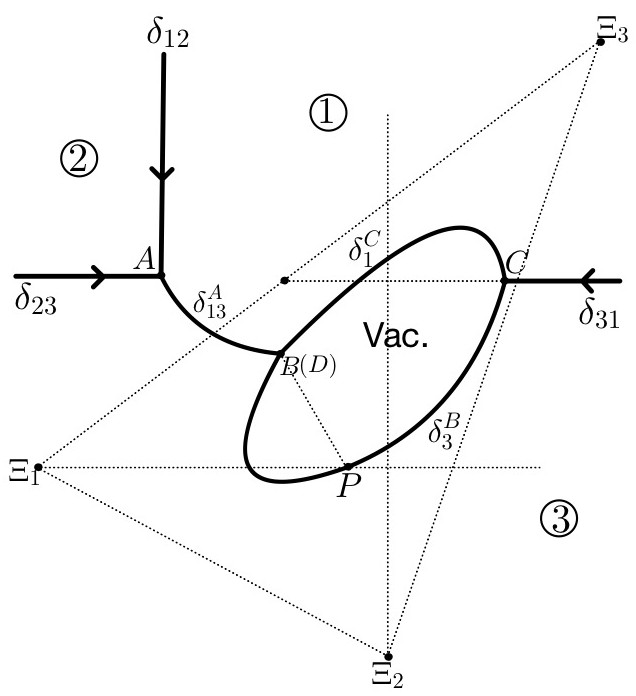}
  \caption{Analytical solution}  
\end{subfigure}%
\begin{subfigure}{.5\textwidth}
  \centering
  \includegraphics[width=.7\linewidth]{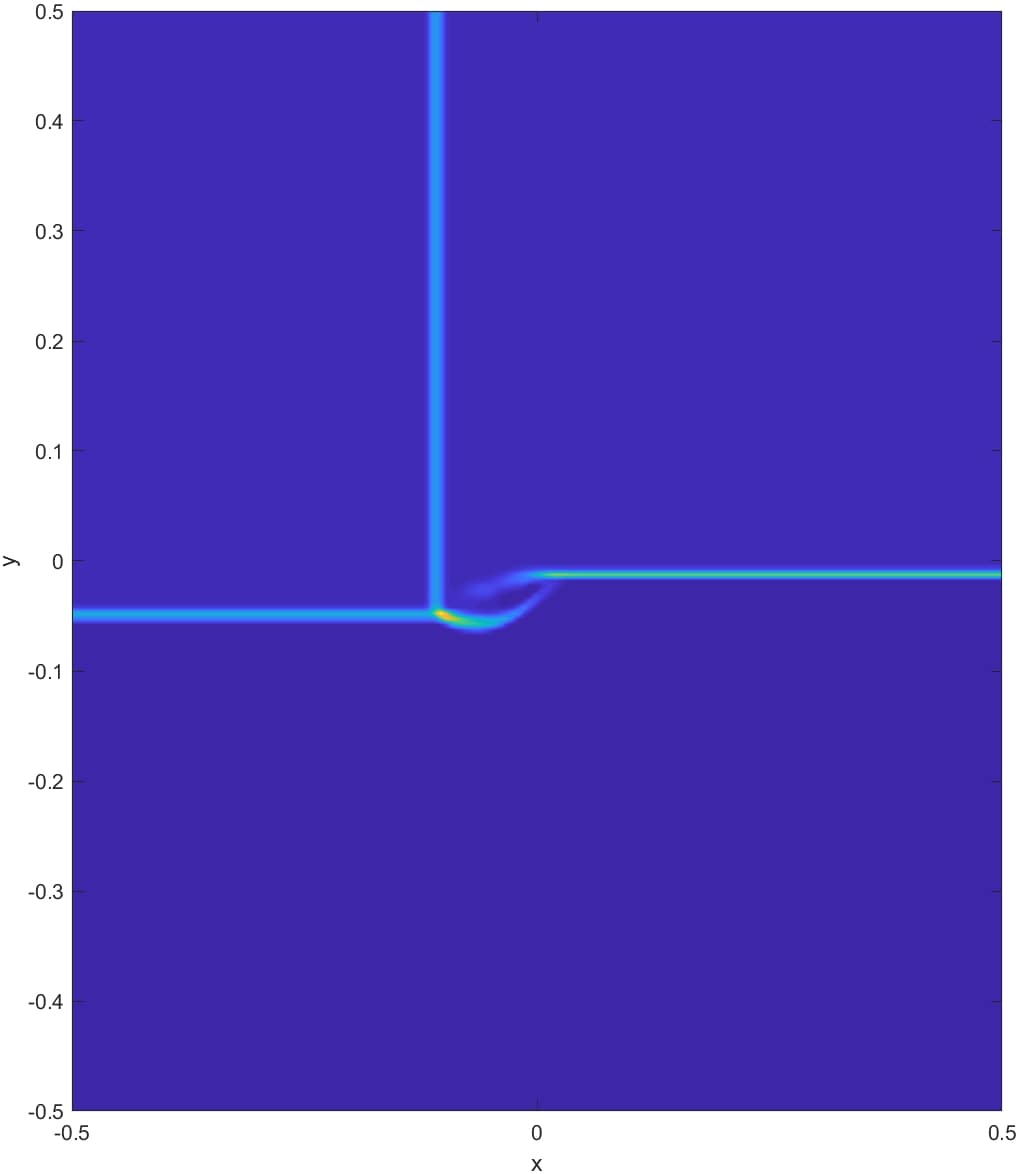}
  \caption{Numerical solution} 
\end{subfigure}
\caption{Solution structure of the Riemann problem (\ref{3}-\ref{4}) via analytical and numerical approach for the case $u_1< u_2 < u_3$, $v_1 < v_{\delta}^{23} < v_3, [\Xi_1, \Xi_3, A] > 0$. The local interaction of $\delta_{12}$ and $\delta_{23}$ generates $\delta_{13}^A$, followed by global interactions of $\delta_{13}^A$ that give rise to a Mach-reflection like pattern with an enclosed vacuum region}
\label{case2}
\end{figure}
Hence, by applying fixed point theorem of Schauder \ref{thm1}, we deduce the existence of a fixed point $(\Xi^*, m_*,n_*, U_\delta^*) \in K$, such that $T(\Xi^*, m_*,n_*, U_\delta^*)=(\Xi^*, m_*,n_*, U_\delta^*)$.

Thus, we constructed the global solution, consisting of a vacuum region enclosed by delta shocks, as depicted in Figure \ref{case1}. To describe this case, the initial data are chosen as
\begin{equation*}
(\rho,u,v)(x,y,0)= \begin{cases}
      (0.1, -0.815, 0.035), & \text{if} \;x>0,y>0, \\
      (0.1, -0.015, -0.825), & \text{if} \;x<0,y>0,\\
      (0.1, 0.85, 0.83), & \text{if} \;y<0,x\in \mathbb{R}.
   \end{cases}
\end{equation*}

It is important to emphasize that the emergence of new delta shocks, which separate the vacuum from the constant states, does not directly arise from the interaction between the external delta shocks. Rather, these external delta shocks undergo a bifurcation process before reaching their respective singular points. This phenomenon exhibits a pattern akin to Mach-reflection, resembling the configuration seen in gas dynamics \cite{jbook}.
 
\subsubsection{Case 2}\label{new}
We begin by assuming that the delta shocks $\delta_{12}$ and $\delta_{23}$ interact at the point $A$. For this configuration, the initial data satisfies the conditions $u_1< u_2 < u_3$, $v_1 < v_{\delta}^{23} < v_3$ and $[\Xi_1, \Xi_3, A] > 0$. The corresponding initial data at this interaction point can be described as follows:
\begin{align}\label{Case2_1}
    \bar{s}=0: \begin{cases}
        \Xi(0)=A=\left(\frac{\sqrt{\rho_1}u_1+\sqrt{\rho_2}u_2}{\sqrt{\rho_1}+\sqrt{\rho_2}},\frac{\sqrt{\rho_2}u_2+\sqrt{\rho_3}u_3}{\sqrt{\rho_3}+\sqrt{\rho_2}} \right),\\
        U_\delta^A(0) = \frac{m_{12}U_\delta^{12}+m_{23}U_\delta^{23}}{m_{12}+m_{23}},\\
        m(0)=m_{12}+m_{23}, \; n(0)=n_{12}+n_{23}\\
        (\rho_L,U_L,H_L)=(\rho_1,U_1,H_1), \; (\rho_R,U_R,H_R)=(\rho_3,U_3,H_3).
    \end{cases}
\end{align}
At this stage, the point $U_2(\Xi_2)$ and the interaction point $A$ lie on opposite sides of the segment $\overline{\Xi_1\Xi}_3$. Consequently, by solving the equation (\ref{pre12}) with the initial data (\ref{Case2_1}), a new delta shock wave, denoted by $\delta_{13}^{A}$, emerges from point $A$ and penetrates the segment $\overline{\Xi_1 \Xi_3}$ in accordance with the behavior described in Lemma \ref{Cased}.  

In this analysis, we focus on the scenario where the trajectory of $\delta_{13}^{A}$ protrudes to the point $\Xi_1$. Under this assumption, two possibilities arise:
\begin{itemize}
    \item $\delta_{13}^{A}$ interacts with $\delta_{31}$, leading to a Riemann problem at the point of interaction. However, this problem does not admit a solution.
    \item Alternatively, the entropy condition for $\delta_{13}^A$ is violated before it reaches to $\delta_{31}$.
\end{itemize}
In either case, the local approach fails to yield a consistent solution. Therefore, a global analysis becomes necessary, ultimately resulting in a Mach-reflection like configuration.

Consider a point $B$ on $\delta_{13}^A$ such that $[\Xi_1,\Xi_2,B]<0$. Suppose a new delta shock, $\delta_3^B$ originates from the point $B$ with initial velocity $U_3^B(0)$ at $B$ satisfying $[\Xi_3,B,U_3^B(0)]>0$ and separates state $\text{\textcircled{\raisebox{-0.9pt}{3}}}$ from the vacuum. Then $\delta_3^B$ must interact with the delta shock $\delta_{31}$  at the point $C$, resulting in the formation of a new delta shock $\delta_1^C$ which connects state  $\text{\textcircled{\raisebox{-0.9pt}{1}}}$ and vacuum. Clearly, the delta shock $\delta_1^C$ interacts either $\delta_{13}^A$ or $\delta_3^B$ at a point $D$ such that $[\Xi_1,\Xi_2,D]<0$. Applying Schauder’s fixed point theorem \ref{thm1}, we establish the existence of a point $B$ on $\delta_{13}^A$ for which $D$  coincides with $B$. This leads us to the construction of a global solution, as shown in Figure \ref{case2}. To illustrate this scenario numerically, we select initial values as
\begin{equation*}
(\rho,u,v)(x,y,0)= \begin{cases}
      (0.5, -0.975, -0.498), & \text{if} \;x>0,y>0, \\
      (0.3, 0.015, -0.925), & \text{if} \;x<0,y>0,\\
      (0.1, 0.945, 0.935), & \text{if} \;y<0,x\in \mathbb{R}.
   \end{cases}
\end{equation*}
\begin{figure}
\begin{subfigure}{.5\textwidth}
  \centering
  \includegraphics[width=.75\linewidth]{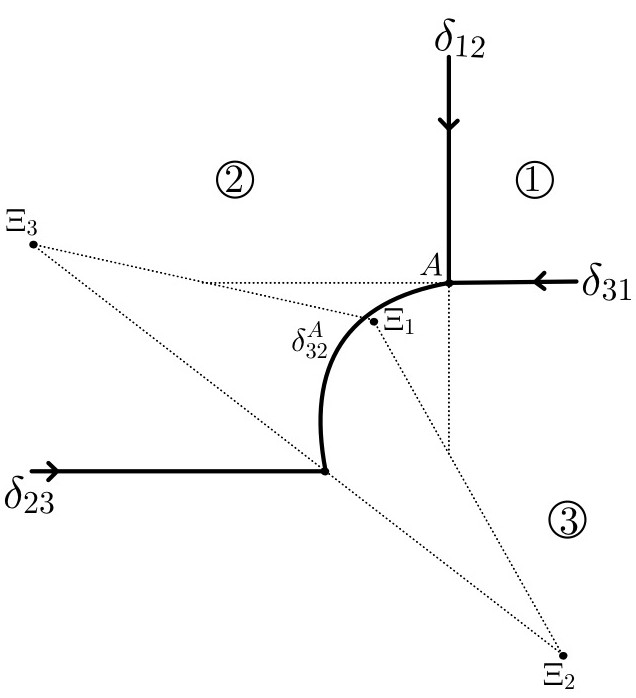}
  \caption{Analytical solution} 
\end{subfigure}%
\begin{subfigure}{.5\textwidth}
  \centering
  \includegraphics[width=.7\linewidth]{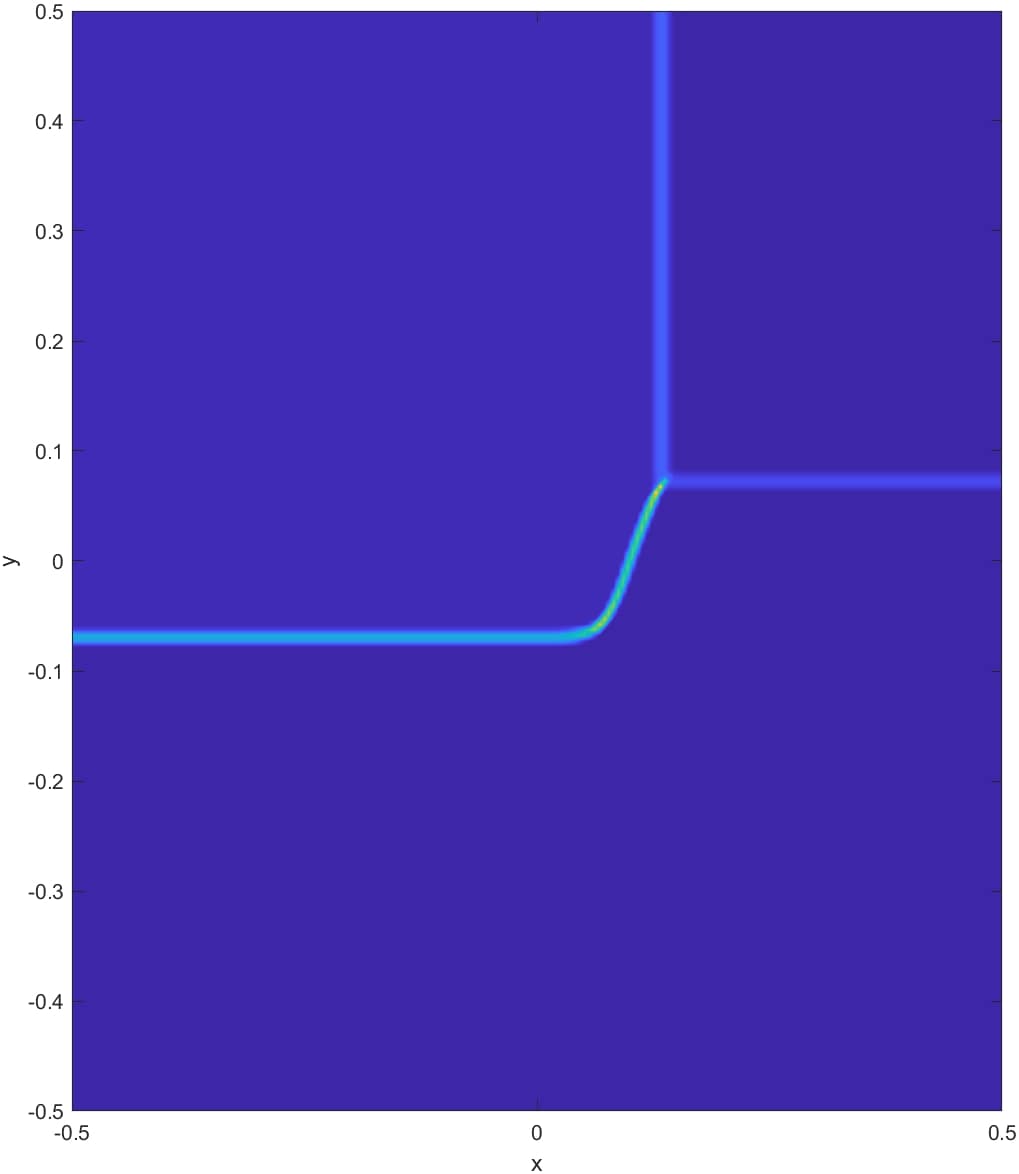}
  \caption{Numerical solution} 
\end{subfigure}
\caption{Analytical and numerical solutions of the Riemann problem (\ref{3})–(\ref{4}) for initial data satisfying $u_3 < u_1 < u_2, [\Xi_1, \Xi_2, \Xi_3] < 0$. The delta shocks $\delta_{12}$ and $\delta_{31}$ interact to form a new delta shock which meets $\delta_{23}$ at its singular point}
\label{case3}
\end{figure}
\subsubsection{Case 3}
In this scenario, the initial configuration is characterized by the inequalities $u_3<u_1<u_2, [\Xi_1,\Xi_2,\Xi_3]<0$. The delta shocks $\delta_{12}$ and $\delta_{31}$ collide at the point $A$ before either reaches its corresponding singularity points. Clearly, $[U_2,U_3,A]<0$.

The associated initial conditions at the interaction point $A$ are given by:
\begin{align}\label{Case3_1}
    \bar{s}=0: \begin{cases}
        \Xi(0)=A=\left(\frac{\sqrt{\rho_1}u_1+\sqrt{\rho_2}u_2}{\sqrt{\rho_1}+\sqrt{\rho_2}},\frac{\sqrt{\rho_1}u_1+\sqrt{\rho_3}u_3}{\sqrt{\rho_3}+\sqrt{\rho_1}} \right),\\
        U_{32}^A(0) = \frac{m_{12}U_\delta^{12}+m_{31}U_\delta^{31}}{m_{12}+m_{31}},\\
        m(0)={m_{12}+m_{31}}, \; n(0)=n_{12}+n_{31}\\
        (\rho_L,U_L,H_L)=(\rho_3,U_3,H_3), \; (\rho_R,U_R,H_R)=(\rho_2,U_2,H_2).
    \end{cases}
\end{align}
This implies that $[U_2,U_3,U_{32}^A]<0$, indicating that both point $A$ and $U_{32}^A(0)$ lie on the same side of the segment $\overline{\Xi_2\Xi}_3$ and moreover, $|[U_2,U_3,A]|>|[U_2,U_3,U_{32}^A]|$. Based on the analysis presented in Lemma \ref{Casec}, it is evident that the delta shock $\delta_{32}^A$ continues smoothly up to the point $\Xi = U_{23}^*$ on $\overline{\Xi_2\Xi}_3$. Simultaneously, the delta shock $\delta_{23}$ originating from infinity, also terminates at $\Xi = U_{23}^*$. The complete solution is illustrated in Figure \ref{case3}. This case is presented using initial conditions
\begin{equation*}
(\rho,u,v)(x,y,0)= \begin{cases}
      (0.2, -0.015, 0.248), & \text{if} \;x>0,y>0, \\
      (0.4, 0.275, -0.925), & \text{if} \;x<0,y>0,\\
      (0.2, -0.985, 0.945), & \text{if} \;y<0,x\in \mathbb{R}.
   \end{cases}
\end{equation*}
In contrast to the Mach-reflection-type configuration, this solution is constructed solely by solving Riemann problems at the interaction points. Such behavior is commonly observed in the interaction of 1-D waves.

\subsection{Cases involving two delta shocks}
Without loss of generality, let us assume that the initial discontinuity connecting states  $\text{\textcircled{\raisebox{-0.9pt}{1}}}$ and $\text{\textcircled{\raisebox{-0.9pt}{2}}}$ generates a contact discontinuity denoted by $J_{12}$, for which the initial data holds $u_1=u_2$. For the remaining two delta shocks, the initial data fulfill the condition $v_1<v_3, v_2<v_3$. Without loss of generality, we further assume $v_2<v_1$. Additionally, we assume, without loss of generality, that
$$u_1=u_2, v_2<v_1<v_3.$$
This setting gives rise to the following three possible subcases:
\begin{figure}
\begin{subfigure}{.5\textwidth}
  \centering
  \includegraphics[width=.75\linewidth]{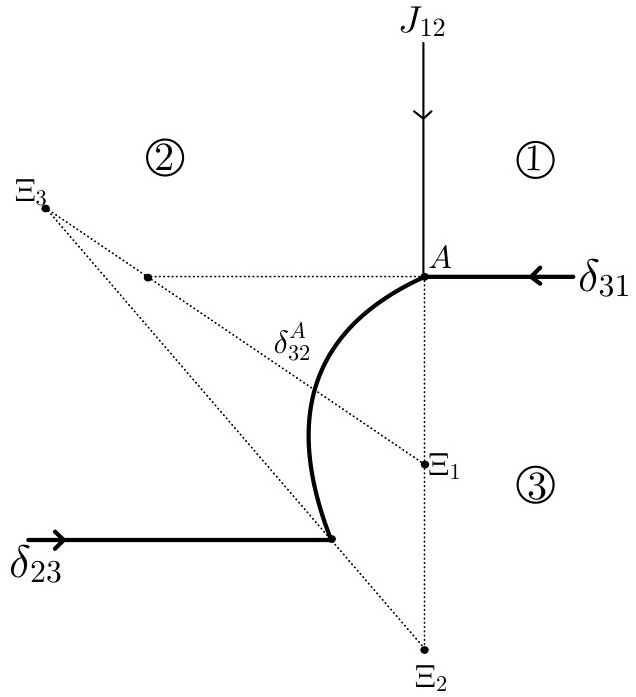}
  \caption{Analytical solution} 
\end{subfigure}%
\begin{subfigure}{.5\textwidth}
  \centering
  \includegraphics[width=.7\linewidth]{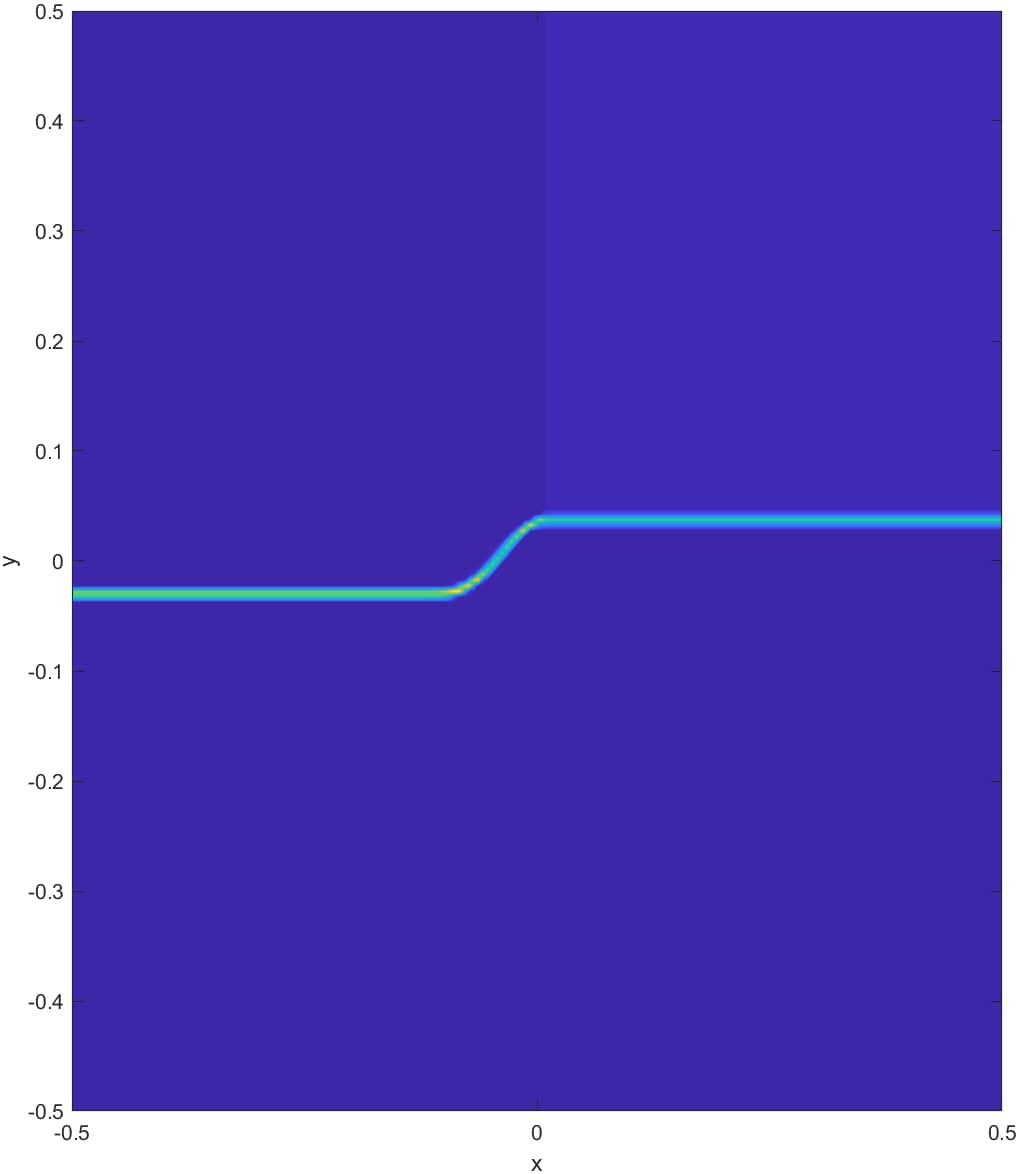}
  \caption{Numerical solution} 
\end{subfigure}
\caption{Analytical and numerical solution of the Riemann problem (\ref{3})–(\ref{4}) for $u_2 = u_1 < u_3, v_2<v_1<v_3$. The contact discontinuity $J_{12}$ interacts with the delta shock $\delta_{31}$ to form a new delta shock, which subsequently meets $\delta_{23}$ at its singular point. Here, the thin line represents the contact discontinuity, while the bold line denotes the delta shock}
\label{case4}
\end{figure}
\subsubsection{Case 4}
This case occurs for the initial data $u_3<u_2=u_1$. The contact discontinuity $J_{12}$ interacts with delta shock $\delta_{31}$ at the point $A$, where initial conditions can be expressed as
\begin{align}\label{Case4_1}
    \bar{s}=0: \begin{cases}
        \Xi(0)=A=(u_1,v_\delta^{31}), \;
        U_{32}^A(0) = U_\delta^{31},\;
        m(0)=m_{31}, \; n(0)=n_{31}\\
        (\rho_L,U_L,H_L)=(\rho_3,U_3,H_3), \; (\rho_R,U_R,H_R)=(\rho_2,U_2,H_2).
    \end{cases}
\end{align}
By analyzing equation (\ref{pre12}) and above initial data (\ref{Case4_1}), we obtain a new delta shock $\delta_{32}^A$, which terminates at $\Xi=U^*_{23}$. The delta shock $\delta_{23}$ meets $\delta_{32}^A$ at $\Xi=U^*_{23}$. Figure \ref{case4} presents the full structure of the resulting solution. For the purpose of this case, the initial values are taken as
\begin{equation*}
(\rho,u,v)(x,y,0)= \begin{cases}
      (0.9, 0.035, -0.24), & \text{if} \;x>0,y>0, \\
      (0.6, 0.035, -0.925), & \text{if} \;x<0,y>0,\\
      (0.3, -0.985, 0.945), & \text{if} \;y<0,x\in \mathbb{R}.
   \end{cases}
\end{equation*}
\subsubsection{Case 5}

\begin{figure}
\begin{subfigure}{.5\textwidth}
  \centering
  \includegraphics[width=.75\linewidth]{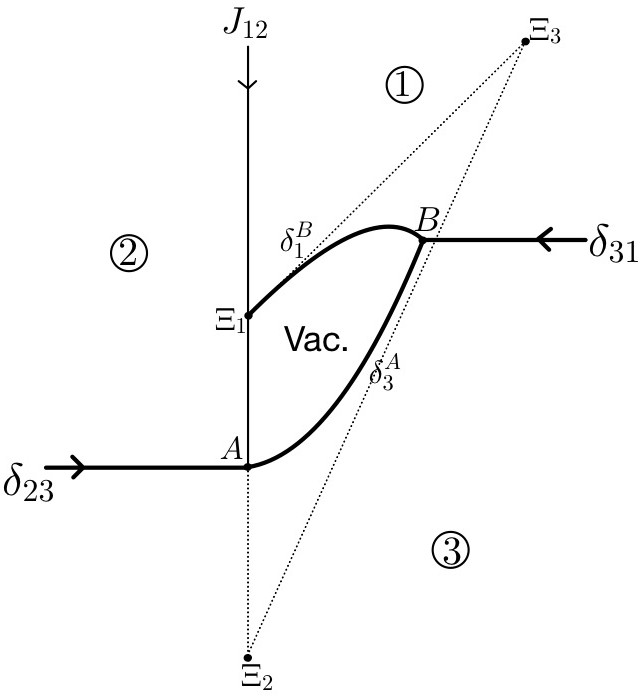}
  \caption{Analytical solution}  
\end{subfigure}%
\begin{subfigure}{.5\textwidth}
  \centering
  \includegraphics[width=.7\linewidth]{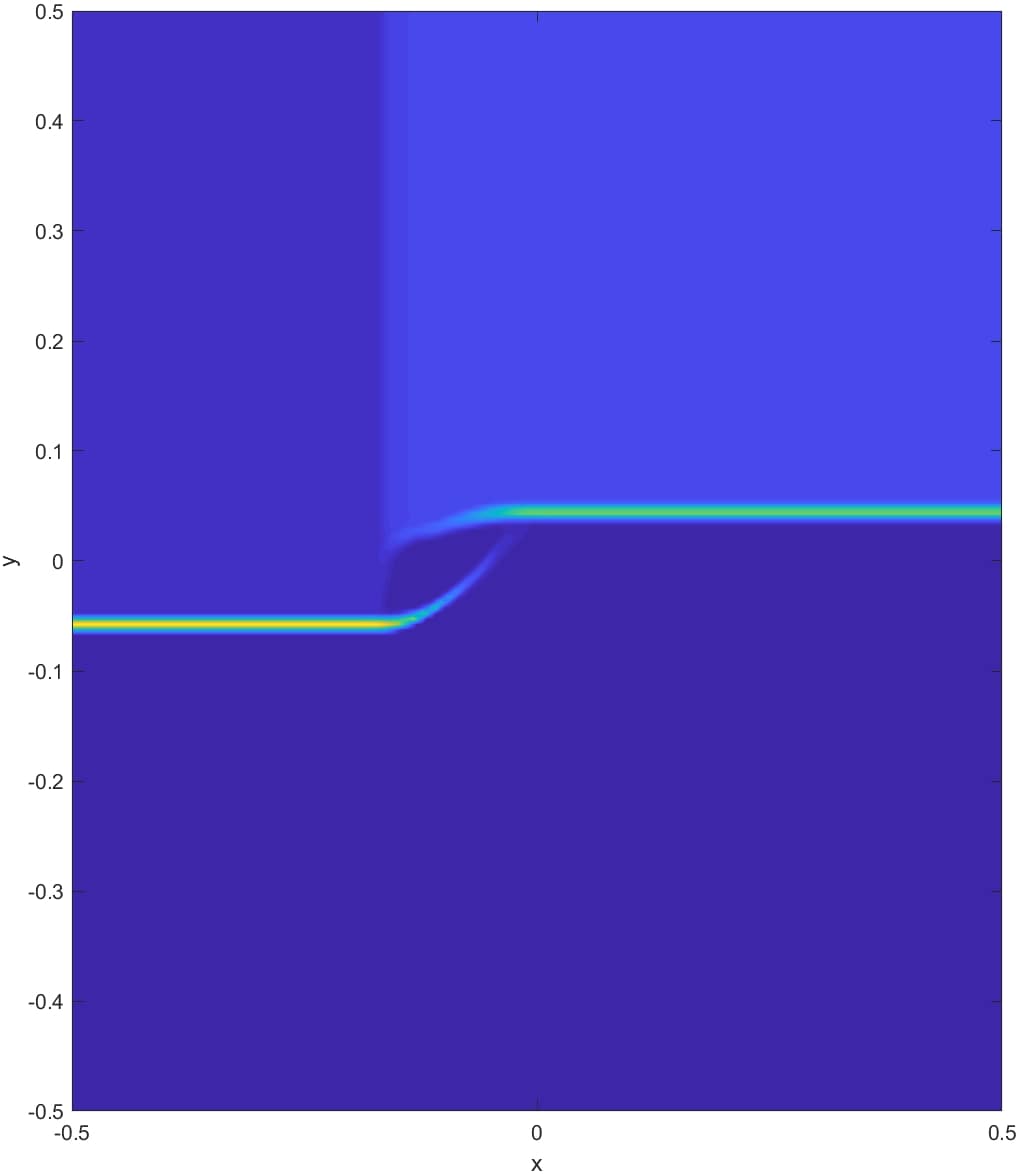}
  \caption{Numerical solution} 
\end{subfigure}
\caption{Analytical and numerical solution of the Riemann problem (\ref{3})–(\ref{4}) for the case $u_2=u_1<u_3, v_2<v_\delta^{23}<v_1$. The contact discontinuity $J_{12}$ interacts with the delta shock $\delta_{23}$, leading to the development of a vacuum region separated by the delta shocks. Here, the thin line represents the contact discontinuity, while the bold line denotes the delta shock}
\label{case5}
\end{figure}

In the present case, the initial data hold $u_2=u_1<u_3, v_2<v_\delta^{23}<v_1$. The delta shock wave $\delta_{23}$ interacts with the contact discontinuity $J_{12}$ at the point $A=(u_1, v_\delta^{23})$, resulting in the emergence of a new delta shock $\delta_{3}^A$, connecting the state  $\text{\textcircled{\raisebox{-0.9pt}{3}}}$ to the vacuum. Following the dynamics described in Lemma \ref{Casea}, this new delta shock $\delta_{3}^A$ further interacts with $\delta_{31}$ at some point $B$, giving rise to another delta shock $\delta_{3}^B$. This delta shock $\delta_{3}^B$ separates the state  $\text{\textcircled{\raisebox{-0.9pt}{1}}}$ from the vacuum and terminates at $\Xi_1$. The final wave pattern is captured in Figure \ref{case5}. The case is constructed numerically using the initial conditions
\begin{equation*}
(\rho,u,v)(x,y,0)= \begin{cases}
      (0.4, -0.615, 0.415), & \text{if} \;x>0,y>0, \\
      (0.8, -0.615, -0.625), & \text{if} \;x<0,y>0,\\
      (0.2, 0.051, 0.415), & \text{if} \;y<0,x\in \mathbb{R}.
   \end{cases}
\end{equation*}
\subsubsection{Case 6}

The initial configuration for this case is defined as $u_1 = u_2<u_3$, $v_2 < v_1 < v_\delta^{23}$. Without loss of generality, we further assume $v_\delta^{23} < v_\delta^{13}$.

Let the interaction point of $J_{12}$ and $\delta_{23}$ is $A$. Solving the Riemann problem at this point yields a delta shock $\delta_{13}^A$. The velocity of $\delta_{13}^A$ at the point $A$ is $U_{13}^A(0)$, which satisfies $[U_1, U_3, U_{13}^A(0)] < 0$, implying that a local solution is not feasible. Hence, similar to Case 2, a global construction of the solution is required.

\begin{figure}
\begin{subfigure}{.5\textwidth}
  \centering
  \includegraphics[width=.75\linewidth]{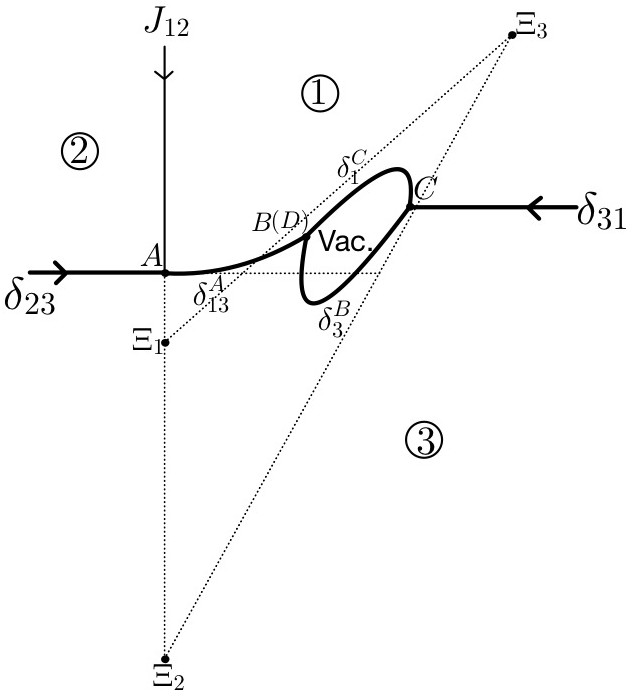}
  \caption{Analytical solution} 
\end{subfigure}%
\begin{subfigure}{.5\textwidth}
  \centering
  \includegraphics[width=.7\linewidth]{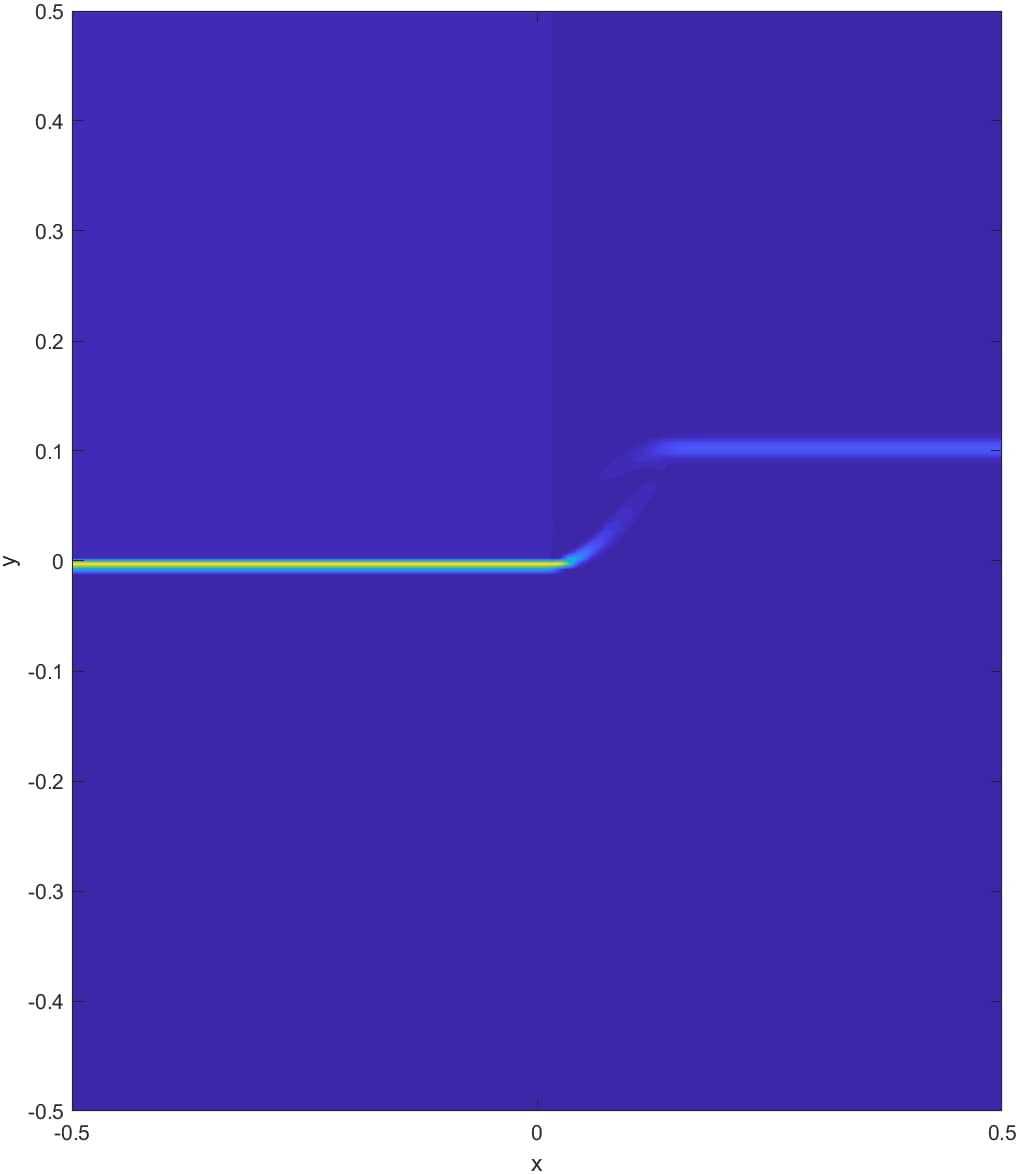}
  \caption{Numerical solution} 
\end{subfigure}
\caption{Solution structure of the Riemann problem (\ref{3}-\ref{4}) via analytical and numerical approach for $u_1 = u_2<u_3$, $v_2 < v_1 < v_\delta^{23}< v_\delta^{13}$. Initially, the delta shock $\delta_{23}$ and contact discontinuity $J_{12}$ interact to produce $\delta_{13}^A$, which subsequently undergoes global interactions, forming a Mach-reflection like configuration with an enclosed vacuum region. Here, the thin line represents the contact discontinuity, while the bold line denotes the delta shocks}
\label{case6}
\end{figure}

Assume that a new delta shock wave $\delta_3^B$ emerges from a point $B$ on $\delta_{13}^A$, where $[U_1, U_3, B] < 0$.
The velocity at this point $B$, denoted by $U_3^B(0)$, satisfies $[\Xi_3, B, U_3^B(0)] > 0$. As a result, $\delta_3^B$ interacts $\delta_{31}$ at some point $C$ from which a new delta shock $\delta_1^C$ originates. This new delta shock $\delta_1^C$ must interact either $\delta_{13}^A$ or $\delta_3^B$ at some point $D$, where $[U_1, U_3, D] < 0$. According to the fixed point theorem of Schauder \ref{thm1}, a point $B$ exists on the delta shock $\delta_{13}^A$ where $D$ and $B$ coincide. Consequently, we arrive at a globally constructed solution similar to that depicted in Figure \ref{case2}. The complete structure is illustrated in Figure \ref{case6}. To simulate this case numerically, the initial conditions are set as
\begin{equation*}
(\rho,u,v)(x,y,0)= \begin{cases}
      (0.2, 0.075, 0.018), & \text{if} \;x>0,y>0, \\
      (0.5, 0.075, -0.625), & \text{if} \;x<0,y>0,\\
      (0.2, 0.985, 0.945), & \text{if} \;y<0,x\in \mathbb{R}.
   \end{cases}
\end{equation*}
\begin{figure}
\begin{subfigure}{.5\textwidth}
  \centering
  \includegraphics[width=.75\linewidth]{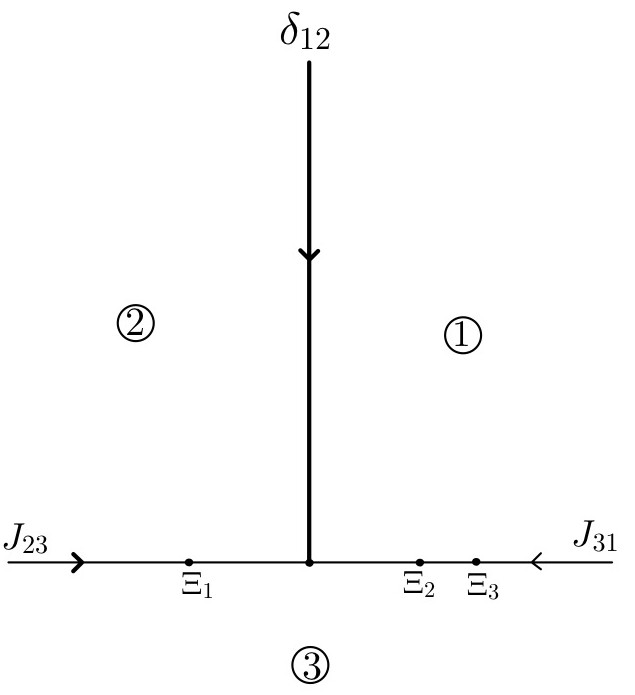}
  \caption{Analytical solution} 
\end{subfigure}%
\begin{subfigure}{.5\textwidth}
  \centering
  \includegraphics[width=.7\linewidth]{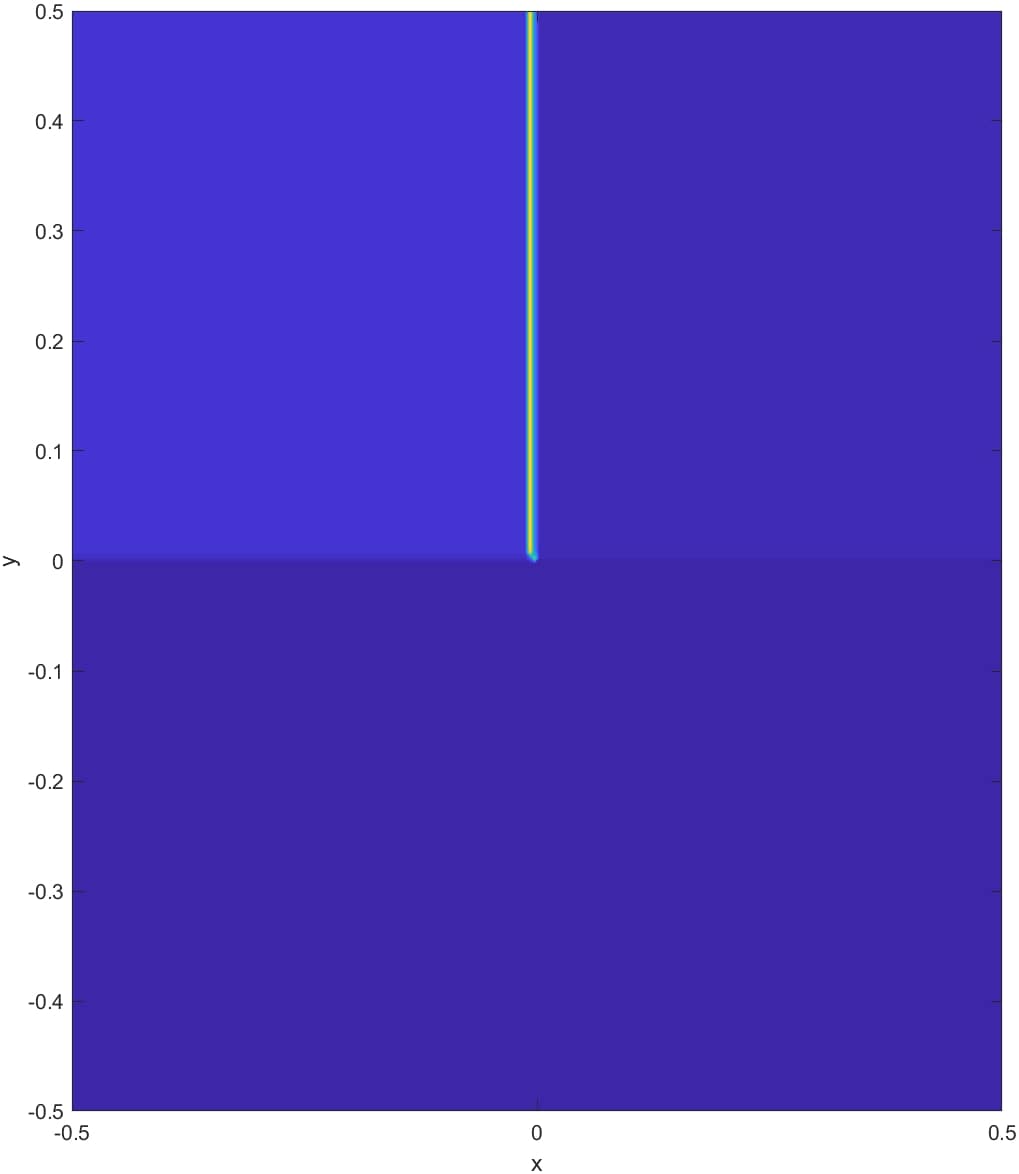}
  \caption{Numerical solution} 
\end{subfigure}
\caption{Solution structure of the Riemann problem (\ref{3}-\ref{4}) via analytical and numerical approach for $v_1 = v_2 = v_3, u_1 < u_2$. The delta shock $\delta_{12}$ interacts with contact discontinuity $J_{23}$, forming a new contact discontinuity which meets the contact discontinuity $J_{31}$. Here, the thin line represents the contact discontinuity, while the bold line denotes the delta shocks}
\label{case7}
\end{figure}
\subsection{Cases involving single delta shock}

Without loss of generality, let the contact discontinuity  $J_{31}$ connects the states  $\text{\textcircled{\raisebox{-0.9pt}{1}}}$ and  $\text{\textcircled{\raisebox{-0.9pt}{3}}}$ with the initial condition $v_3=v_1$. 
In this configuration, three distinct solution structures may arise: 

\subsubsection{Case 7}
Without loss of generality, consider that the states $\text{\textcircled{\raisebox{-0.9pt}{1}}}$ and  $\text{\textcircled{\raisebox{-0.9pt}{2}}}$ to be connected by the delta shock $\delta_{12}$. In this configuration, the initial values obey $v_1 = v_2 = v_3, u_1 < u_2$. The contact discontinuity $J_{23}$ meets the delta shock $\delta_{12}$ at the point $U_{12}^*$, resulting in the emergence of a new contact discontinuity $J_{13}$ which subsequently interacts the contact discontinuity $J_{31}$ at the point $\Xi = \Xi_3$. In addition, $J_{31}$ and $J_{13}$ lies on the same line. The final wave pattern is captured in Figure \ref{case7}. To describe this case, initial values are prescribed as
\begin{equation*}
(\rho,u,v)(x,y,0)= \begin{cases}
      (0.2, -0.475, 0.018), & \text{if} \;x>0,y>0, \\
      (0.3, 0.275, 0.018), & \text{if} \;x<0,y>0,\\
      (0.1, 0.435, 0.018), & \text{if} \;y<0,x\in \mathbb{R}.
   \end{cases}
\end{equation*}
\begin{figure}
\begin{subfigure}{.5\textwidth}
  \centering
  \includegraphics[width=.75\linewidth]{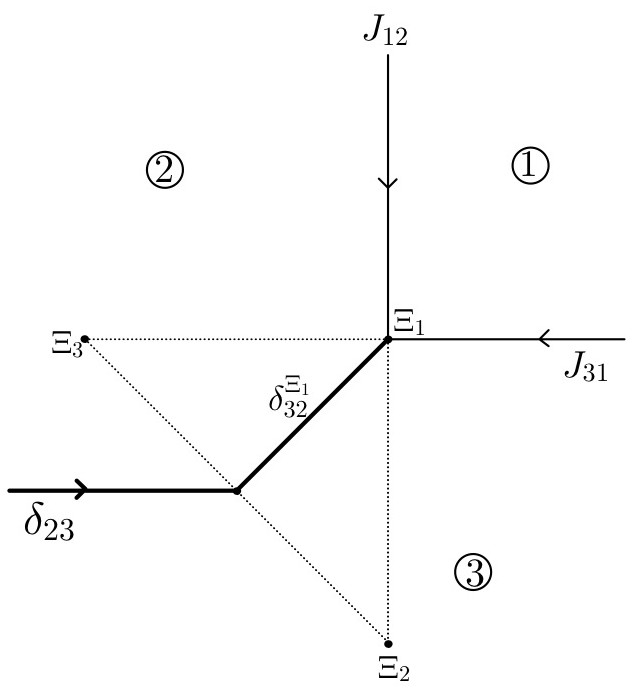}
  \caption{Analytical solution}  
\end{subfigure}%
\begin{subfigure}{.5\textwidth}
  \centering
  \includegraphics[width=.7\linewidth]{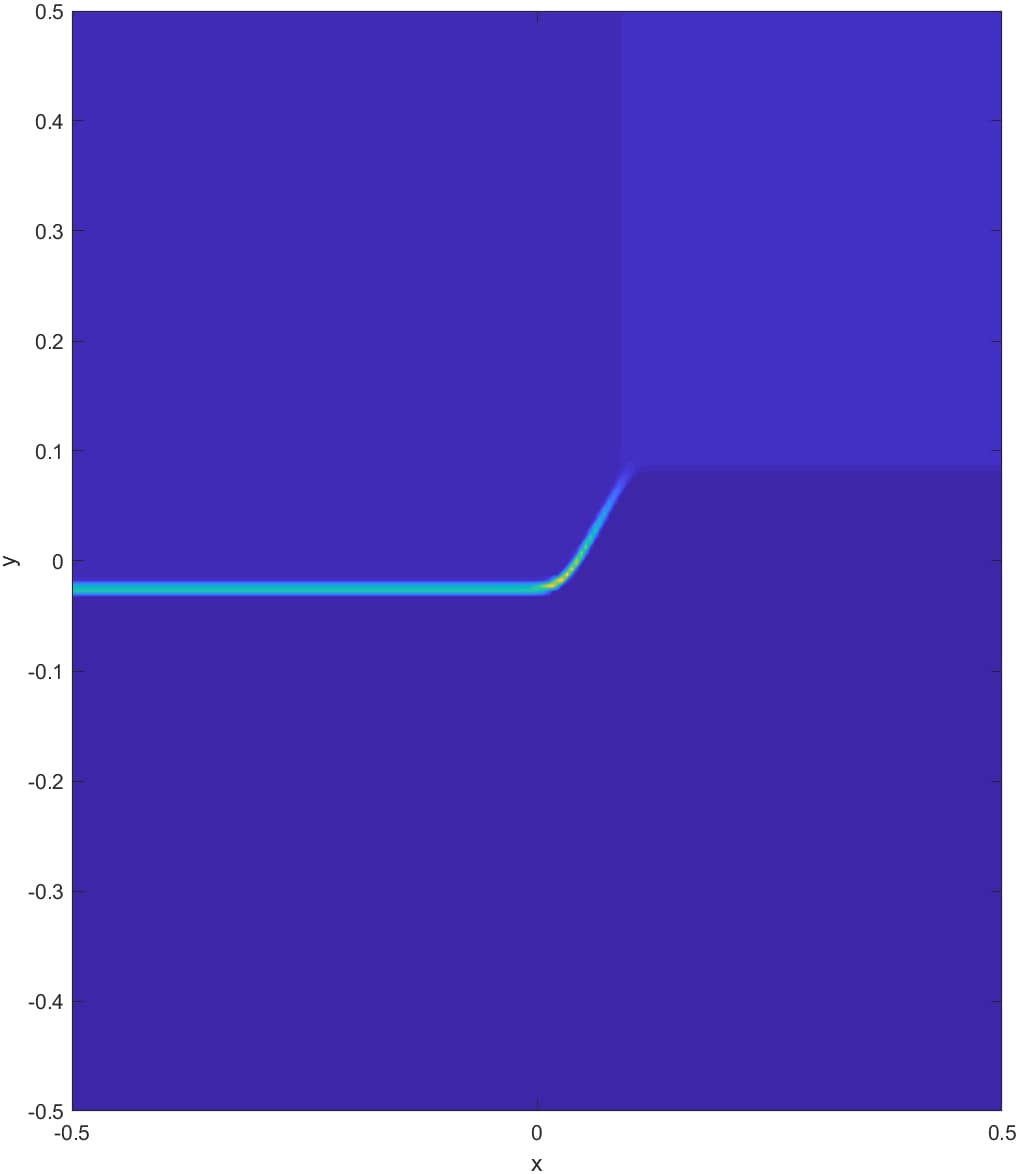}
  \caption{Numerical solution} 
\end{subfigure}
\caption{Solution structure of the Riemann problem (\ref{3}-\ref{4}) via analytical and numerical approaches for $u_3<u_1=u_2, v_2<v_1=v_3$. The contact discontinuities $J_{12}$ and $J_{31}$ interacts forming a delta shock which meets the delta shock $\delta_{23}$ at $ U_{23}^*$. Here, the thin line represents the contact discontinuity, while the bold line denotes the delta shocks}
\label{case8}
\end{figure}
\subsubsection{Case 8}\label{new1}
For this case, assume that the delta shock $\delta_{23}$ connects the states  $\text{\textcircled{\raisebox{-0.9pt}{2}}}$ and  $\text{\textcircled{\raisebox{-0.9pt}{3}}}$. We first consider the scenario where $u_3<u_1=u_2$. Hence, the initial condition holds the relation $u_3<u_1=u_2, v_2<v_1=v_3$. The contact discontinuities $J_{12}$ and $J_{31}$ interact at the point $\Xi_1$. Solving equation (\ref{pre12}) at $\Xi_1$, yields a new delta shock $\delta_{32}^{\Xi_1}$ which then interacts with the delta shock $\delta_{23}$ at $U^*_{23}$. This interaction leads to the structure displayed in Figure \ref{case8}. For the analysis of this case, the initial data is assumed to be
\begin{equation*}
(\rho,u,v)(x,y,0)= \begin{cases}
      (0.9, 0.475, 0.418), & \text{if} \;x>0,y>0, \\
      (0.2, 0.475, -0.438), & \text{if} \;x<0,y>0,\\
      (0.6, -0.435, 0.418), & \text{if} \;y<0,x\in \mathbb{R}.
   \end{cases}
\end{equation*}

\subsubsection{Case 9}
\begin{figure}
\begin{subfigure}{.5\textwidth}
  \centering
  \includegraphics[width=.75\linewidth]{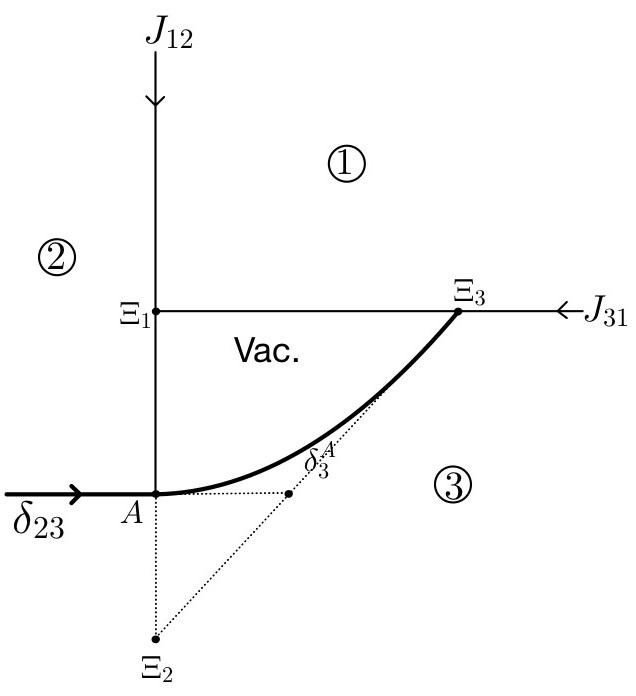}
  \caption{Analytical solution}  
\end{subfigure}%
\begin{subfigure}{.5\textwidth}
  \centering
  \includegraphics[width=.7\linewidth]{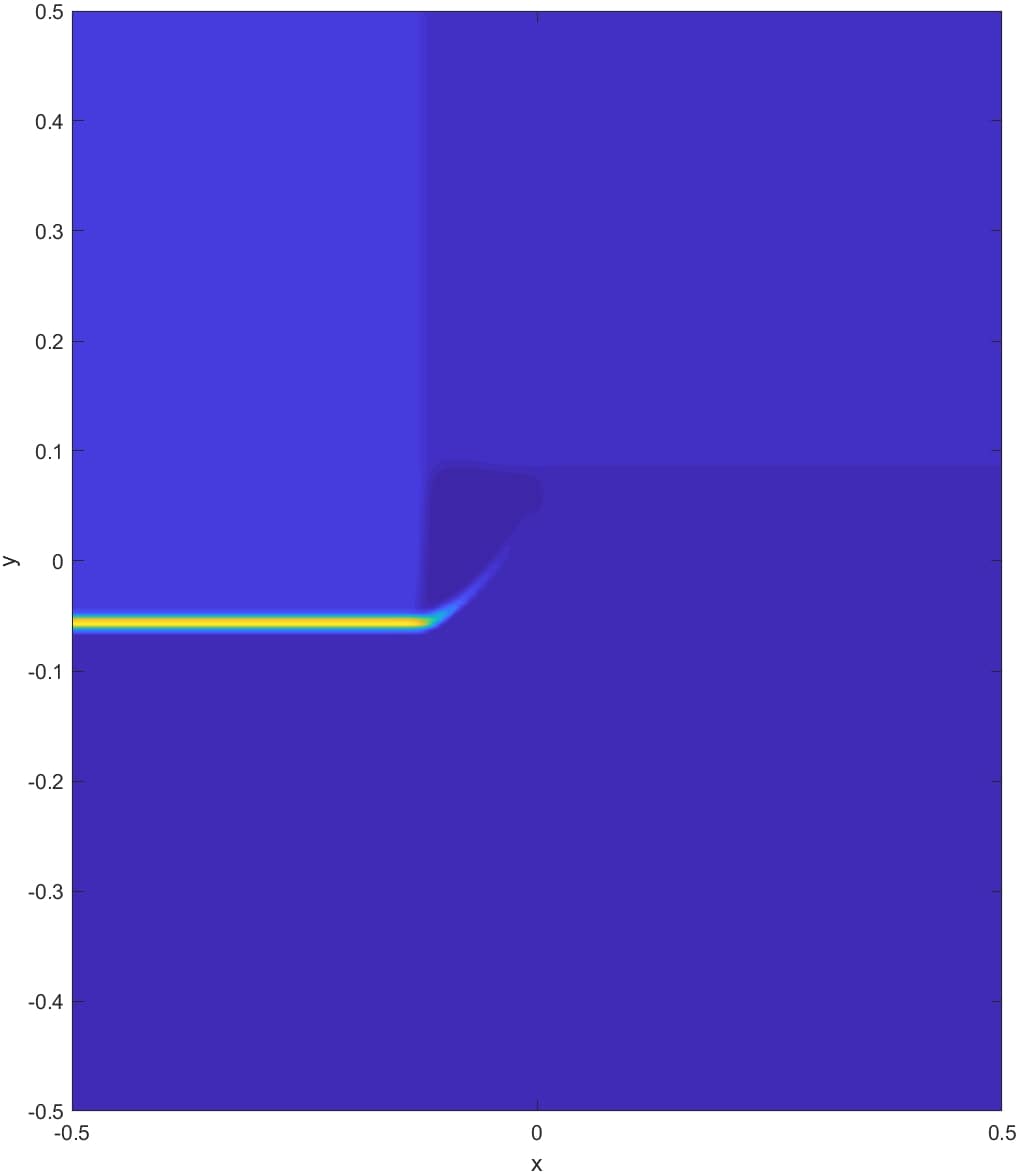}
  \caption{Numerical solution} 
\end{subfigure}
\caption{Analytical and numerical representation of the solution of Riemann problem (\ref{3}-\ref{4}) for $u_1=u_2<u_3, v_2<v_1=v_3$. A vacuum region emerges in the solution due to the interaction between the contact discontinuity $J_{12}$ and the delta shock $\delta_{23}$. The thin line represents the contact discontinuity, while the bold line indicates the delta shock} 
\label{case9}
\end{figure}

Continuing from Case 8, we now consider the reverse configuration where the initial data satisfy $u_1=u_2<u_3, v_2<v_1=v_3$. In this setting, the delta shock $\delta_{23}$ interacts the contact discontinuity at a point $A$, triggering the occurrence of a new delta shock $\delta_3^A$ that connects state  $\text{\textcircled{\raisebox{-0.9pt}{3}}}$ to the vacuum. This delta shock ultimately terminates at $\Xi_3$. Figure \ref{case9} demonstrates the resulting configuration of the solution. The initial data is specified to perform the numerical simulation for this case
\begin{equation*}
(\rho,u,v)(x,y,0)= \begin{cases}
      (0.4, -0.615, 0.415), & \text{if} \;x>0,y>0, \\
      (0.8, -0.615, -0.625), & \text{if} \;x<0,y>0,\\
      (0.2, 0.051, 0.415), & \text{if} \;y<0,x\in \mathbb{R}.
   \end{cases}
\end{equation*}

\section{Conclusions}

We thoroughly analyzed the solution of the two-dimensional Riemann problem (\ref{3}-\ref{4}) under the assumption that each jump of the discontinuity projects exactly one wave, resulting in topologically distinct nine solutions along with their corresponding criteria. Interactions involving contact discontinuities and delta shocks are examined by employing the method of characteristic analysis. In some solution structures, a vacuum region appears. It is important to note that the occurrence of a delta shock wave, characterized by a Dirac delta function in both density and internal energy, results from the overlapping of linearly degenerate characteristic lines. Moreover Mach-reflection like configurations are observed, where a single delta shock splits into two distinct delta shocks, forming a triple wave point where all three such waves interact. The analytical solutions are further validated numerically using a second-order, semidiscrete central-upwind scheme on a uniform mesh of $200 \times 200$ cells and numerical results are obtained at time $T=0.20$. It is observed that the obtained numerical solutions have good agreement with the corresponding analytical solutions. In future, we focus on formulating a robust second and higher-order Godunov-type solver for the two-dimensional pressureless gas dynamics system with generalized Riemann initial data and investigate its convergence.

\section*{Acknowledgements}
\textit{The first author acknowledges the financial support provided by the Indian Institute of Technology Kharagpur. The second author (TRS) expresses his gratitude towards SERB, DST, Government of India (Ref. No. CRG/2022/006297) for its financial support through the core research grant.}

\end{document}